\documentclass[11pt,a4paper]{article}
\usepackage[T1]{fontenc}
\usepackage{lmodern}
\usepackage[margin=27mm]{geometry}
\usepackage{amsmath,amssymb,amsthm,mathtools}
\usepackage{microtype}
\usepackage[hidelinks]{hyperref}

\allowdisplaybreaks[1]
\numberwithin{equation}{section}
\newtheorem{theorem}{Theorem}[section]
\newtheorem{corollary}[theorem]{Corollary}
\newtheorem{proposition}[theorem]{Proposition}
\newtheorem{lemma}[theorem]{Lemma}
\theoremstyle{remark}

\DeclareMathOperator{\diver}{div}
\DeclareMathOperator{\supp}{supp}

\title{A sharp convergence criterion for a Yamabe-type flow}

\author{%
Hongyu Chen\thanks{School of Mathematics, Sichuan University,
Chengdu 610064, China.
Email: \texttt{hongyu.chern@gmail.com}.}
\qquad
Li Ma\thanks{Corresponding author.
School of Mathematics and Physics,
University of Science and Technology Beijing,
Beijing 100083, China.
Email: \texttt{nuslma@126.com}.}
}
\date{October 2026}

\begin{document}
\maketitle
\vspace{-2em}

\begin{abstract}
Let $(M^n,g)$ be a closed connected Riemannian manifold and $\Psi$ a
smooth potential with positive generalized Yamabe constant. For $n\ge4$,
we prove that the normalized Yamabe-type flow converges smoothly from
every positive smooth initial value if and only if $\Psi\le R_g$.
 In dimension three, either
$\Psi\le R_g$ or positivity of the Green function's regular part at
every pole ensures convergence at arbitrary initial energy. On the
unit round $\mathbb{RP}^3$, all trajectories converge for constant
potentials $6<\Psi<8$, which exceed the scalar curvature. In dimension four, at any prescribed point where $\Psi>R_g$, we
construct a concentrating trajectory with initial energy arbitrarily
close to the spherical value from above.
\end{abstract}

\section{Introduction}\label{introduction}
  In this work, we study the asymptotic behavior of Yamabe type flow, which is a global flow as showed in \cite{Ma23}. This flow is a generalization of the Yamabe flow studied by S.Brendle in \cite{Bre05} and \cite{Bre07}, where the author have used the positive mass theorem \cite{LP87}.

Let $(M^n,g)$ be a closed connected smooth Riemannian manifold of
dimension $n\ge3$, with scalar curvature $R_g$. We study convergence
of the normalized flow obtained by replacing $R_g$ in the conformal
Laplacian by a smooth potential $\Psi$. We use the convention that
$\Delta_g v\le0$ at a local maximum and put
$$
p=\frac{n+2}{n-2},\qquad c_n=\frac{4(n-1)}{n-2},\qquad
L_\Psi=-c_n\Delta_g+\Psi.
$$
The associated energy and generalized Yamabe constant are
\begin{equation}\label{energy}
E_\Psi(v)=
\frac{\displaystyle\int_M(c_n|\nabla v|_g^2+\Psi v^2)\,dv_g}
{\displaystyle\left(\int_M |v|^{p+1}\,dv_g\right)^{2/(p+1)}},
\qquad
Y_\Psi(M)=\inf_{0\ne v\in H^1(M)} E_\Psi(v).
\end{equation}
The normalized Yamabe-type flow introduced by Ma \cite{Ma23} is
\begin{equation}\label{flow}
u_t=-u^{1-p}L_\Psi u+r(t)u,\qquad
r(t)=\int_M(c_n|\nabla u|_g^2+\Psi u^2)\,dv_g,\qquad u(0)=u_0,
\end{equation}
with $u_0>0$ smooth and $\int_M u_0^{p+1}\,dv_g=1$. The flow
preserves this mass and decreases $E_\Psi$. Write $Y=Y(S^n)$
for the Yamabe constant of the round sphere, which is the energy
of a mass-one Euclidean bubble. We ask whether the comparison
$\Psi\le R_g$ determines convergence from every positive smooth
initial value in the regime $Y_\Psi(M)>0$.

One of our mains is below, which is about the asymptotic behavior of the global Yamabe type flow for dimension bigger than three.
\begin{theorem}\label{main-iff}
Suppose $n\ge4$ and $Y_\Psi(M)>0$. The condition $\Psi\le R_g$
on $M$ is equivalent to convergence in every $C^m(M)$ of the solutions
of \eqref{flow} from all positive smooth initial values of critical
mass one.

If $\Psi(y)>R_g(y)$, then for every $\eta>0$ there is a positive smooth
initial value of critical mass one with $Y<E_\Psi(u_0)<Y+\eta$
such that
$$
 r(t)\longrightarrow Y,\qquad
 u(t)^{p+1}\,dv_g\rightharpoonup\delta_y,\qquad
 \|u(t)\|_{L^\infty(M)}\longrightarrow\infty .
$$
\end{theorem}

In dimension four, the leading energy excess of a cutoff bubble
of scale $\epsilon$ at $y$ is a positive multiple of
$(\Psi(y)-R_g(y))\epsilon^2\log(1/\epsilon)$. Where $\Psi>R_g$ this
dominates the squared $L^{4/3}$ norm of the bubble's residual, which
is $O(\epsilon^2)$. The scale derivative bounds the excess by the
square root of the dissipation. This controls the total variation
of the critical-mass barycenter, whose uniform tail bound allows us
to prescribe its limit by Brouwer's theorem.

In dimension three, the leading energy correction for a projected
bubble is determined by the regular part of the Green function.
When $Y_\Psi(M)>0$, normalize $G_\Psi$ by
$L_{\Psi,x}G_\Psi(x,y)=\delta_y$ and set
\begin{equation}\label{three-mass}
 A_\Psi(y)=\lim_{x\to y}\left(G_\Psi(x,y)
                  -\frac{1}{32\pi\,d_g(x,y)}\right).
\end{equation}
The limit exists and is continuous in $y$ by
Section~\ref{projective-section}.

Our result in dimension three can be stated as follows.
\begin{theorem}\label{main-three}
Suppose $n=3$ and $Y_\Psi(M)>0$. If $A_\Psi(y)>0$ at every point,
then every positive smooth initial value of critical mass one gives
a solution of \eqref{flow} which converges in $C^m(M)$ for every
integer $m\ge0$. The same conclusion holds if $\Psi\le R_g$.
\end{theorem}

\begin{corollary}\label{main-counterexample}
Let $M=\mathbb{RP}^3$ with the metric induced by the unit round
three-sphere. For any constant potential $\Psi\in(6,8)$ the flow
\eqref{flow} converges smoothly from every positive smooth initial
value of critical mass one.
\end{corollary}

Here $R_g=6$, so these are examples of convergence with $\Psi>R_g$
everywhere. The positive solutions in the three statements exist
uniquely for all $t\ge0$, and the limits in the convergence assertions
satisfy
\begin{equation}\label{stationary}
 L_\Psi u_\infty=r_\infty u_\infty^p,\qquad
 u_\infty>0,\qquad
 \int_Mu_\infty^{p+1}\,dv_g=1,\qquad
 r_\infty=\lim_{t\to\infty}r(t)>0.
\end{equation}

At arbitrary initial energy, the convergence proof has to allow
several bubbles and a possible positive stationary part. We derive a
gradient inequality from Green projections whose squared
$L^{2n/(n+2)}$ residual is of lower order than the energy gain.
For $n\ge4$ and $\Psi\le R_g$, the resolvent identity contributes at
least $c\epsilon^{n-2}\int_M(R_g-\Psi)\,dv_g$ to the gain, with $c>0$
independent of the scale and center. Thus strict inequality somewhere
gives a gain even at centers where $\Psi=R_g$ on a neighborhood.

At $\Psi=R_g$, \eqref{flow} is the classical Yamabe flow after a
constant time rescaling. Its convergence follows from the results
of Ye and Brendle \cite{Ye94,Bre05,Bre07}, with the higher-dimensional
mass condition described in \cite[Section~5]{BreSurvey}.
Ma obtained a positive stationary subsequential limit when
$0<Y_\Psi(M)<Y$ and $E_\Psi(u_0)\le Y$ \cite[Theorem~2]{Ma23}.
For $n\ge5$, the concentration assertion also follows from
Kim--Musso \cite[Theorem~1.1]{KM24} with perturbation $\Psi-R_g$
after normalizing the critical mass and rescaling time.

Integrals and norms on $M$ use $dv_g$. We set $a=(n-2)/2$ and
$s=2n/(n+2)=(p+1)/p$. Uniform estimates are for the fixed manifold
and potential and the indicated scales and centers.

The remaining part of this paper is organized as follows. We establish the evolution estimates and
bubble decomposition in Section~\ref{evolution-section}. We prove in Section~\ref{source-section} the
convergence criterion used in Sections~\ref{all-dimensional-bubbles}
and in \ref{projective-section} for the potential comparison and the
three-dimensional Green function, respectively.
In Section~\ref{four-section} we construct the concentrating trajectories.

\section{Evolution and global compactness}\label{evolution-section}

The assumption $Y_\Psi(M)>0$ makes the quadratic form of $L_\Psi$ an equivalent Hilbert norm on $H^1(M)$. We write
$$
\langle v,w\rangle_\Psi=\int_M(c_n\langle\nabla v,\nabla w\rangle_g+\Psi vw).
$$
To check coercivity, \eqref{energy} and H\"older's inequality give $\|v\|_2^2\le C\|v\|_{p+1}^2\le C\langle v,v\rangle_\Psi$. The identity defining the quadratic form then controls $\|\nabla v\|_2^2$ even when $\Psi$ changes sign. The reverse norm bound follows from boundedness of $\Psi$.
We write $\|v\|_\Psi^2=\langle v,v\rangle_\Psi$.
Coercivity also gives a positive Green kernel $G_\Psi$ with
$c\,d_g(x,y)^{2-n}\le G_\Psi(x,y)\le C\,d_g(x,y)^{2-n}$
off the diagonal. To see positivity, divide by a positive first
eigenfunction and apply the maximum principle. The local
fundamental singularity gives the bounds near the diagonal,
and positivity and compactness give them on its complement.

\begin{proposition}\label{evolution}
Assume $Y_\Psi(M)>0$. The smooth positive solution of \eqref{flow} exists for all $t\ge0$, preserves $\int_Mu^{p+1}=1$, and satisfies
\begin{equation}\label{dissipation}
r'(t)=-2\int_Mu^{p-1}u_t^2,
\qquad Y_\Psi(M)\le r(t)\le r(0).
\end{equation}
Its limit $r_\infty>0$ exists, and
\begin{equation}\label{residual-decay}
\lim_{t\to\infty}\int_Mu^{p-1}u_t^2=0,
\qquad
\lim_{t\to\infty}\|L_\Psi u-r_\infty u^p\|_{2n/(n+2)}=0.
\end{equation}
\end{proposition}

\begin{proof}
Local existence holds for this quasilinear parabolic equation on
the open set of positive functions. Multiplication by $u^p$ gives
$$
\frac1{p+1}\frac d{dt}\int_Mu^{p+1}
=r\left(\int_Mu^{p+1}-1\right).
$$
The initial normalization and local uniqueness for this scalar equation near the positive value one give $\int_Mu^{p+1}=1$. Differentiating $r$ and using this identity proves \eqref{dissipation}.

Ma's proof of global existence uses the pseudo-scalar curvature
\cite[Theorem~1 and Section~2]{Ma23}. With $g(t)=u^{4/(n-2)}g$, write
$$
R_\Psi=u^{-p}L_\Psi u,
\qquad u_t=u(r-R_\Psi).
$$
The product rule for $L_\Psi(uv)$ and the conformal formula for the Laplacian give
\begin{equation}\label{curvature-evolution}
\partial_tR_\Psi
=c_n\Delta_{g(t)}R_\Psi+(p-1)R_\Psi(R_\Psi-r).
\end{equation}
At a negative minimum the reaction term is nonnegative, since $r>0$. The maximum principle yields
\begin{equation}\label{curvature-lower}
R_\Psi\ge-C,\qquad u_t\le(r(0)+C)u.
\end{equation}
Thus $u$ has a finite upper bound on every finite time interval. If $u\le A$ on such an interval, choose a positive constant $K>\sup\Psi+CA^{p-1}$. Then $(-c_n\Delta_g+K)u\ge0$. Its positive Green kernel has a positive lower bound on $M\times M$: near the diagonal the kernel tends to infinity, and off the diagonal the assertion follows from positivity and compactness. Consequently,
$$
u(x)\ge c_KK\int_Mu\ge c_KK A^{-p}.
$$
The equation is thus uniformly parabolic on any finite
time interval. The change of dependent variable and time
used in Lemma \ref{smoothing} removes the nonlocal reaction,
and the positive bounds give the local continuation
estimates. Local uniqueness extends to the global solution.

Since $dv_{g(t)}=u^{p+1}dv_g$, differentiation of this measure and \eqref{curvature-evolution} give
\begin{align}\label{variance-evolution}
\frac d{dt}\int_M(R_\Psi-r)^2\,dv_{g(t)}
={}&-2c_n\int_M|\nabla(R_\Psi-r)|_{g(t)}^2\,dv_{g(t)}\notag\\
&+2(p-1)r\int_M(R_\Psi-r)^2\,dv_{g(t)}\notag\\
&+(p-3)\int_M(R_\Psi-r)^3\,dv_{g(t)}.
\end{align}
Here the term containing $r'$ vanishes because $\int_M(R_\Psi-r)\,dv_{g(t)}=0$. The quantity on the left before differentiation is $\int_Mu^{p-1}u_t^2$, whose time integral is finite by \eqref{dissipation}.

If $n\ge4$, then $p-3\le0$. The lower bound $R_\Psi-r\ge-C$ bounds the last term in \eqref{variance-evolution} from above by $C\int_M(R_\Psi-r)^2\,dv_{g(t)}$. Hence its derivative is bounded above by the same quadratic integral times a constant.

For $n=3$ one has $p=5$. The exact quadratic-form identity
$$
\int_M(uh)L_\Psi(uh)
=\int_M\bigl(c_n|\nabla h|_{g(t)}^2+R_\Psi h^2\bigr)\,dv_{g(t)}
$$
and the definition of $Y_\Psi$ give, with $h=R_\Psi-r$,
\begin{align*}
Y_\Psi\|R_\Psi-r\|_{L^6(g(t))}^2
\le{}&c_n\int_M|\nabla(R_\Psi-r)|_{g(t)}^2\,dv_{g(t)}\\
&+r\int_M(R_\Psi-r)^2\,dv_{g(t)}
+\int_M(R_\Psi-r)^3\,dv_{g(t)}.
\end{align*}
Interpolation on the probability space $(M,dv_{g(t)})$ yields
$$
\int_M|R_\Psi-r|^3\,dv_{g(t)}
\le\left(\int_M(R_\Psi-r)^2\,dv_{g(t)}\right)^{3/4}
\|R_\Psi-r\|_{L^6(g(t))}^{3/2}.
$$
Young's inequality first absorbs half the $L^6$ term in the preceding Sobolev inequality. Substituting that bound into this interpolation inequality, and applying Young's inequality to the gradient term, gives, for every fixed $\varepsilon>0$,
\begin{align*}
\int_M|R_\Psi-r|^3\,dv_{g(t)}
\le{}&\varepsilon c_n\int_M|\nabla(R_\Psi-r)|_{g(t)}^2\,dv_{g(t)}\\
&+C_\varepsilon\left[
\left(\int_M(R_\Psi-r)^2\,dv_{g(t)}\right)^{3/2}
+\left(\int_M(R_\Psi-r)^2\,dv_{g(t)}\right)^3\right].
\end{align*}
It follows that the derivative of the quadratic integral is at most a constant times the sum of its first, $3/2$ and third powers. While that integral is at most one, its derivative is therefore at most a constant times itself. It must eventually stay below one: every late upcrossing from $1/2$ to one would, by integration of this differential inequality, contribute a fixed positive amount to its time integral. This contradicts integrability.

In every dimension the quadratic integral is nonnegative and integrable, and, for all sufficiently large times, its derivative is bounded above by a constant multiple of itself. It therefore tends to zero. Otherwise a value at least $\varepsilon$ forces values at least $\varepsilon/2$ on a fixed preceding time interval, contrary to the vanishing tail integral. Finally, H\"older's inequality and $\int_Mu^{p+1}=1$ give
\begin{equation}\label{residual-dissipation}
\|L_\Psi u-r_\infty u^p\|_{2n/(n+2)}
\le\left(\int_Mu^{p-1}u_t^2\right)^{1/2}+r-r_\infty.
\end{equation}
The exponent in the second H\"older factor is $(p-1)n/2=p+1$. This proves \eqref{residual-decay}.
\end{proof}

The mass-one Euclidean Sobolev extremals \cite{AubSob76,T76} are
\begin{equation}\label{euclidean-bubble}
 U_{\epsilon,y}(x)=K
 \left(\frac{\epsilon}{\epsilon^2+|x-y|^2}\right)^a,\qquad
 K=\left(\frac{4n(n-1)}{Y}\right)^{(n-2)/4}.
\end{equation}
They satisfy $\int_{\mathbb R^n} U_{\epsilon,y}^{p+1}=1$ and
$-c_n\Delta U_{\epsilon,y}=Y U_{\epsilon,y}^p$. We write
$U_\epsilon=U_{\epsilon,0}$, also as a function of the radial
variable.

\begin{proposition}\label{decomposition}
Let $u_j>0$ be smooth, $\int_Mu_j^{p+1}=1$, and suppose that $u_j$ is bounded in $H^1(M)$ and
$$
L_\Psi u_j-\lambda u_j^p\longrightarrow0
\quad\text{in }L^{2n/(n+2)}(M),\qquad \lambda>0.
$$
After taking a subsequence, there are an integer $k\ge0$, a smooth nonnegative solution $u_*$ of $L_\Psi u_*=\lambda u_*^p$, centers $y_{i,j}$ and scales $\epsilon_{i,j}\to0$, such that $u_j$ is the sum of $u_*$, $k$ cut-off local standard bubbles with coefficient $(Y/\lambda)^{(n-2)/4}$, and a remainder tending to zero in $H^1$. The parameters satisfy
\begin{equation}\label{orthogonal-scales}
\frac{\epsilon_{i,j}}{\epsilon_{\ell,j}}
+\frac{\epsilon_{\ell,j}}{\epsilon_{i,j}}
+\frac{d_g(y_{i,j},y_{\ell,j})^2}{\epsilon_{i,j}\epsilon_{\ell,j}}
\longrightarrow\infty\quad(i\ne\ell),
\end{equation}
and
\begin{equation}\label{mass-splitting}
\int_Mu_*^{p+1}+k\left(\frac{Y}\lambda\right)^{n/2}=1.
\end{equation}
The function $u_*$ is either zero or strictly positive.
\end{proposition}

\begin{proof}
The nonnegative rescaled sequence $(\lambda/c_n)^{1/(p-1)}u_j$
is Palais--Smale for the critical functional with operator
$-\Delta_g+\Psi/c_n$ and nonlinear coefficient one. Its residual
tends to zero in $H^{-1}$ by Sobolev duality, and testing the original
residual against $u_j$ gives $E_\Psi(u_j)\to\lambda$ and hence
convergence of the rescaled action. We use Mazumdar's closed-manifold
extension \cite{Maz17} of Struwe's decomposition \cite{Str84}.
Theorem~1.2 and Proposition~5.1 of the
\href{https://arxiv.org/abs/1603.07953v1}{preprint version}
give the decomposition and parameter separation. The interior
argument uses normal coordinates and the Riemannian density on
any closed manifold.
The positive full-space profiles are the bubbles of \cite{CGS89}.
Using a fixed local cutoff changes each by a vanishing Dirichlet
tail. The Brezis--Lieb lemma \cite[Theorem~1]{BL83} gives the
mass splitting.
Scaling the coefficient of each profile back to $\lambda$
gives its amplitude $(Y/\lambda)^{1/(p-1)}$ and its
mass $(Y/\lambda)^{n/2}$. A nonzero weak background
is first made bounded by the local critical elliptic
iteration. In that iteration the $L^{n/2}$ norm of
$\lambda u_*^{p-1}$ on sufficiently small balls absorbs
the potential term at any fixed finite set of exponents.
Elliptic regularity and the strong maximum principle
give a smooth positive background.
\end{proof}

\section{Convergence from projected sources}\label{source-section}

We follow Brendle's finite-dimensional reduction \cite{Bre05},
using Green projections of nonnegative sources for the bubbles.
Related gradient estimates for the Euclidean Dirichlet problem
appear in \cite[Theorem~4.1 and Section~5]{JXZ26}.

The main result in this section is the following one.
\begin{proposition}\label{source-convergence}
Let $(M^n,g)$ be closed, connected and smooth, $n\ge3$, and let
$\Psi\in C^\infty(M)$ satisfy $Y_\Psi(M)>0$.
Suppose there is a smooth family of nonnegative smooth sources
$f_{\epsilon,y}$, $y\in M$, $0<\epsilon<\epsilon_0$, with
$B_{\epsilon,y}=YL_\Psi^{-1}f_{\epsilon,y}$, satisfying the
following estimates uniformly in $y$:
\begin{gather}
\|f_{\epsilon,y}\|_s=1,\qquad
f_{\epsilon,y}(x)\le
C\left(\frac{\epsilon}{\epsilon^2+d_g(x,y)^2}\right)^{(n+2)/2},
\qquad
\int_{B_\epsilon(y)}f_{\epsilon,y}\ge c\epsilon^a,
\label{family-profile}\\
D_\Psi(\epsilon,y):=\int_Mf_{\epsilon,y}L_\Psi^{-1}f_{\epsilon,y}-Y^{-1}
>0,\qquad D_\Psi(\epsilon,y)\longrightarrow0,
\label{family-gain}\\
\|L_\Psi B_{\epsilon,y}-YB_{\epsilon,y}^p\|_s^2
=o(D_\Psi(\epsilon,y)).
\label{family-error}
\end{gather}
Assume also that $B_{\epsilon,y}$ and its first two normalized
parameter derivatives, formed from $\epsilon\partial_\epsilon$
and $\epsilon\nabla_y$, differ by $o(1)$ in $H^1(M)$ from the
corresponding mass-one Euclidean bubble with a fixed local cutoff
and its parameter derivatives.
Then every positive smooth mass-one initial value for
\eqref{flow} gives a global solution which converges smoothly
to a positive solution of \eqref{stationary}.
\end{proposition}

The parameter statement is understood in any of a finite set of
smooth local frames; a change of frame has uniformly bounded
coefficients and gives the same assertion. To prove Proposition \ref{source-convergence}, we need to set up new estimates in following subsections.

\subsection{Interactions between bubbles}\label{bubble-energy}

We assume the hypotheses of Proposition \ref{source-convergence}
throughout this section.
The Green-kernel estimate in Section \ref{evolution-section} and \eqref{family-profile} give
$$
0<B_{\epsilon,y}(x)\le
C\left(\frac{\epsilon}{\epsilon^2+d_g(x,y)^2}\right)^a.
$$
Indeed, for $r=d_g(x,y)>2\epsilon$, split the source integral
into $B_{r/2}(y)$, $B_{r/2}(x)$ and the remaining dyadic annuli.
Their bounds are, respectively,
$C\epsilon^ar^{-2a}$, $C\epsilon^{a+2}r^{-2a-2}$ and
$C\epsilon^ar^{-2a}$. For $r\le2\epsilon$, rescaling gives
$C\epsilon^{-a}$.

For a fixed finite family put $B_i=YL_\Psi^{-1}f_{\epsilon_i,y_i}$ and
$$
\varepsilon_{i\ell}=
\left(\frac{\epsilon_i}{\epsilon_\ell}
+\frac{\epsilon_\ell}{\epsilon_i}
+\frac{d_g(y_i,y_\ell)^2}{\epsilon_i\epsilon_\ell}\right)^{-a}.
$$
These numbers tend to zero under \eqref{orthogonal-scales}.

\begin{lemma}\label{pair-estimates}
For two Euclidean bubbles with $\varepsilon_{i\ell}\to0$,
\begin{gather}
\int_{\mathbb R^n}U_i^pU_\ell
=\int_{\mathbb R^n}U_iU_\ell^p\asymp\varepsilon_{i\ell},\notag\\
\int_{\mathbb R^n}\min(U_i^pU_\ell,U_iU_\ell^p)
\le C\varepsilon_{i\ell}^{n/(n-2)}.\label{pair-overlap}
\end{gather}
The actual nonlinear interaction satisfies
\begin{equation}\label{pair-nonlinear}
\|(U_i+U_\ell)^p-U_i^p-U_\ell^p\|_{2n/(n+2)}
\le C\begin{cases}
\varepsilon_{i\ell},&3\le n\le5,\\
\varepsilon_{i\ell}(1+|\log\varepsilon_{i\ell}|)^{2/3},&n=6,\\
\varepsilon_{i\ell}^{(n+2)/(2(n-2))},&n\ge7.
\end{cases}
\end{equation}
For their projected sources and bubbles,
\begin{align}\label{source-pair}
\int_Mf_iL_\Psi^{-1}f_\ell&\ge c\varepsilon_{i\ell},\notag\\
\int_M\left(\sum_i b_if_i\right)^{2n/(n+2)}
&=\sum_i b_i^{2n/(n+2)}
+O\left(\sum_{i<\ell}\varepsilon_{i\ell}^{n/(n-2)}\right),
\end{align}
uniformly for amplitudes in a fixed compact subset of $(0,\infty)$.
\end{lemma}

\begin{proof}
Translate and dilate the first bubble to $U_{1,0}$. Its conformal compactification is a volume-one round sphere. A simultaneous conformal transformation sends the ratio of the second bubble to the first to
$$
b_\mu(x)=\left(\frac{\mu(1+|x|^2)}{\mu^2+|x|^2}\right)^{(n-2)/2},
\qquad 0<\mu\le1,
$$
where $\mu^{(n-2)/2}\asymp\varepsilon_{i\ell}$. The ratio of two conformal bubble factors on the sphere has a denominator affine in the spherical coordinates, and rotation on the sphere gives the displayed $b_\mu$, with
$$
\mu+\mu^{-1}
=\frac{\epsilon_i^2+\epsilon_\ell^2+|y_i-y_\ell|^2}{\epsilon_i\epsilon_\ell}.
$$ Integrals of total degree $p+1$, including the powers defining the norm in \eqref{pair-nonlinear}, are invariant under this common transformation.

On the sphere, $\int b_\mu\asymp\mu^{(n-2)/2}$. The lower bound comes from a fixed annulus, and the upper bound follows by integrating $\mu^{(n-2)/2}(\mu^2+r^2)^{-(n-2)/2}$ for $r\le1$ and retaining the sphere-volume weight for $r\ge1$. The equality of the two original pair integrals also follows by integration by parts in the bubble equations. Since $b_\mu\ge1$ precisely for $r\le\sqrt\mu$, splitting there gives
\begin{align*}
\int\min(b_\mu,b_\mu^p)
\le{}&C\mu^{(n-2)/2}\int_0^{\sqrt\mu}
\frac{r^{n-1}}{(\mu^2+r^2)^{(n-2)/2}}\,dr\\
&+C\mu^{(n+2)/2}\int_{\sqrt\mu}^1r^{-3}\,dr
+C\mu^{(n+2)/2}
\le C\mu^{n/2}.
\end{align*}
This proves \eqref{pair-overlap}.

For $n<6$, $p>2$ and the nonlinear difference is at most $C(U_i^{p-1}U_\ell+U_iU_\ell^{p-1})$. The sphere integrals of the two terms give norms at most $C\mu^{(n-2)/2}$. The required radial conditions are $(n-2)2n/(n+2)<n$ and $4\cdot2n/(n+2)>n$, both valid precisely when $n<6$.

For $n\ge6$, $1<p\le2$ and
$$
0\le(A+B)^p-A^p-B^p\le C\min(A^{p-1}B,AB^{p-1}).
$$
On the sphere the last expression becomes $C\min(b_\mu,b_\mu^{p-1})$. Its $2n/(n+2)$ power is bounded on the two annuli by
\begin{align*}
&C\mu^{4n/(n+2)}r^{-8n/(n+2)}
&&\text{if }\mu<r<\sqrt\mu,\\
&C\mu^{n(n-2)/(n+2)}r^{-2n(n-2)/(n+2)}
&&\text{if }\sqrt\mu<r<1.
\end{align*} Integrating with $r^{n-1}dr$ gives a logarithm when $n=6$ and $O(\mu^{n/2})$ when $n>6$. The core $r<\mu$ and the sphere region $r>1$ are no larger. Taking the $2n/(n+2)$ root proves \eqref{pair-nonlinear}. In every case its right side squared is $o(\varepsilon_{i\ell})$.

For the actual sources, \eqref{family-profile} gives mass at least
$c\epsilon_i^a$ on the core ball of radius $\epsilon_i$. On the product of two such balls the kernel is at least
$c(d_g(y_i,y_\ell)+\epsilon_i+\epsilon_\ell)^{-2a}$. Hence
$$
\int_Mf_iL_\Psi^{-1}f_\ell
\ge
\frac{c(\epsilon_i\epsilon_\ell)^a}
{(d_g(y_i,y_\ell)+\epsilon_i+\epsilon_\ell)^{2a}}
\ge c\varepsilon_{i\ell}.
$$
This argument permits both nested scales and colliding centers.

For $1<s<2$ the excess of the $s$-power of a finite positive sum
is bounded by a fixed multiple of
$$
\sum_{i<\ell}\min(f_i^{s-1}f_\ell,f_if_\ell^{s-1}).
$$
In a common coordinate neighborhood the profile bounds are
$f_i\le C U_{\epsilon_i,y_i}^p$ and
$B_i\le C U_{\epsilon_i,y_i}$, with the Euclidean bubbles expressed
in that one chart. Distances and volumes are uniformly comparable.
Since $p(s-1)=1$, \eqref{pair-overlap} gives the source excess.
For distinct limiting centers, on a neighborhood of $y_i$ the other
source is at most $C\epsilon_\ell^{a+2}$. Compare it with a Euclidean
bubble centered a fixed distance outside the smaller chart.
The same integral estimate applies. On the complement of these
neighborhoods the excess is bounded by
$$
C\min(\epsilon_i^a\epsilon_\ell^{a+2},
       \epsilon_i^{a+2}\epsilon_\ell^a)
\le C(\epsilon_i\epsilon_\ell)^{a+1},
$$
which is comparable to $\varepsilon_{i\ell}^{n/(n-2)}$.
This proves \eqref{source-pair} for the globally supported sources.

The nonlinear power excess is increasing in each positive argument.
The same local comparison therefore proves \eqref{pair-nonlinear}
with $B_i,B_\ell$ in place of $U_i,U_\ell$. Outside all center
neighborhoods the global tails give the same bound. More explicitly,
for $n>6$, $1<p<2$ and the exterior term is bounded by
$C\min(\epsilon_i^2\epsilon_\ell^a,
\epsilon_i^a\epsilon_\ell^2)
\le C(\epsilon_i\epsilon_\ell)^{(a+2)/2}$.
For $n=6$ the exterior estimate is $C(\epsilon_i\epsilon_\ell)^2$;
for $3\le n\le5$ it is $C(\epsilon_i\epsilon_\ell)^a$.
These are no larger than the claimed pair bounds for separated centers.
\end{proof}

\begin{lemma}\label{sum-errors}
For fixed positive bounded amplitudes, the squared $L^{2n/(n+2)}$ norm of
$$
\left(\sum_i b_iB_i\right)^p-\sum_i b_i^pY^{-1}L_\Psi B_i
$$
is
\begin{equation}\label{pure-error}
o\left(\sum_i\left(\int_Mf_iL_\Psi^{-1}f_i-\frac1{Y}\right)
+\sum_{i<\ell}\varepsilon_{i\ell}\right).
\end{equation}
If $u_z$ ranges over a fixed smooth family with uniform positive upper and lower bounds, the squared norm of
$$
\left(u_z+\sum_i b_iB_i\right)^p-u_z^p
-\sum_i b_i^pY^{-1}L_\Psi B_i
$$
is
\begin{equation}\label{mixed-error}
o\left(\sum_i\epsilon_i^{(n-2)/2}
+\sum_i\left(\int_Mf_iL_\Psi^{-1}f_i-\frac1{Y}\right)
+\sum_{i<\ell}\varepsilon_{i\ell}\right).
\end{equation}
One also has
\begin{equation}\label{background-weight-error}
\left\|\left(u_z+\sum_i b_iB_i\right)^{p-1}-u_z^{p-1}\right\|_{2n/(n+2)}^2
=o\left(\sum_i\epsilon_i^{(n-2)/2}\right).
\end{equation}
\end{lemma}

\begin{proof}
The individual equation errors have squared norms
$o(D_\Psi(\epsilon_i,y_i))$ by \eqref{family-error}.
The nonlinear excess of a fixed finite positive sum is bounded by
a fixed multiple of the sum of its pair excesses.
Lemma \ref{pair-estimates} gives squared pair errors
$o(\varepsilon_{i\ell})$ in every dimension.
The triangle inequality, the fixed number of terms and these
two comparisons prove \eqref{pure-error}.

For a positive bounded background, the additional nonlinear error
and the Hessian-weight change are bounded by finite sums of
$C(B_i+B_i^{p-1})$ for $n<6$, and
$C\min(B_i,B_i^{p-1})$ for $n\ge6$.
Splitting at $U_i=1$, whose radius is comparable to
$\sqrt{\epsilon_i}$, gives their $L^s$ norms at most
$$
C\begin{cases}
\epsilon_i^a,&3\le n\le5,\\
\epsilon_i^2(1+|\log\epsilon_i|)^{2/3},&n=6,\\
\epsilon_i^{(n+2)/4},&n>6.
\end{cases}
$$
For $n<6$, integrate $U_i^s$ and $U_i^{(p-1)s}$ individually;
the inequalities $(n-2)s<n<4s$ give the displayed order.
For $n\ge6$, on the inner annulus use $U_i^{p-1}$ and on the
outer annulus use $U_i$. The radial integrals are
$$
\epsilon_i^{2s}\int_{\epsilon_i}^{\sqrt{\epsilon_i}}
r^{n-1-4s}\,dr,\qquad
\epsilon_i^{as}\int_{\sqrt{\epsilon_i}}^{D}
r^{n-1-2as}\,dr.
$$
Here $D$ is a fixed coordinate radius, and the omitted core and
global tail terms have no larger order.
They give a logarithm at $n=6$ and the power
$\epsilon_i^{n/2}$ before taking the $s$-root for $n>6$.
Each resulting squared norm is $o(\epsilon_i^a)$.
This proves \eqref{mixed-error} and \eqref{background-weight-error}.

For later use, the dual power overlap of $f_i$ with a smooth
positive bounded source is at most $C\int\min(U_i,U_i^p)$.
Indeed,
$$
\int_{r<\sqrt{\epsilon_i}}U_i
\le C\epsilon_i^a\int_0^{\sqrt{\epsilon_i}}
\frac{r^{n-1}}{(\epsilon_i^2+r^2)^a}\,dr
\le C\epsilon_i^{n/2},
$$
and
$$
\int_{r>\sqrt{\epsilon_i}}U_i^p
\le C\epsilon_i^{a+2}\int_{\sqrt{\epsilon_i}}^D r^{-3}\,dr
+C\epsilon_i^{a+2}
\le C\epsilon_i^{n/2}.
$$
The negative kernel contribution is at least
$\int_Mu_zf_i\ge c\epsilon_i^a$, by the core mass and the
positive lower bound for $u_z$. The overlap is consequently
smaller than this contribution.
\end{proof}

Fix $\lambda>0$ and consider the action
$$
I_\lambda(v)=\frac12\langle v,v\rangle_\Psi-\frac\lambda{p+1}\int_M|v|^{p+1}.
$$
For every nonnegative source $F$, Young's inequality with conjugate exponents $p+1$ and $(p+1)/p$ gives
\begin{equation}\label{dual-action}
I_\lambda(L_\Psi^{-1}F)
\le-\frac12\int_MF L_\Psi^{-1}F
+\frac{p}{p+1}\lambda^{-1/p}\int_MF^{(p+1)/p}.
\end{equation}
The negative sign of the kernel term will be retained.

\subsection{Modulation and the local gradient inequality}\label{modulation-section}
 The below is about the finite dimensional reduction.
\begin{lemma}\label{analytic-background}
Let $u_*>0$ solve $L_\Psi u_*=\lambda u_*^p$. There is a finite-dimensional subspace $P\subset H^1(M)$ of smooth functions, an orthonormal basis $\phi_1,\ldots,\phi_N$ for $\langle\cdot,\cdot\rangle_\Psi$, and a real analytic family
$$
u_z=u_*+\sum_{\ell=1}^Nz_\ell\phi_\ell+\zeta(z),
\qquad \zeta(z)\in P^\perp,\qquad \zeta(0)=D\zeta(0)=0,
$$
defined near $z=0$, such that $u_z>0$, $L_\Psi u_z>0$, and
\begin{equation}\label{reduced-equation}
(1-\Pi)\bigl(u_z-\lambda L_\Psi^{-1}(u_z^p)\bigr)=0.
\end{equation}
Here $\Pi$ is the orthogonal projection onto $P$. The Hilbert Hessian of $I_\lambda$ is strictly positive on $P^\perp$ at $u_*$. For some $\gamma\in(0,1]$,
\begin{equation}\label{finite-dimensional-gradient}
|I_\lambda(u_z)-I_\lambda(u_*)|
\le C|\nabla_z I_\lambda(u_z)|^{1+\gamma}.
\end{equation}
For $w\perp_\Psi P$, one has the exact cancellation
\begin{equation}\label{tangent-cancellation}
\int_M(L_\Psi w-p\lambda u_z^{p-1}w)\,\partial_{z_\ell}u_z=0.
\end{equation}
\end{lemma}

\begin{proof}
The operator $h\mapsto L_\Psi^{-1}(u_*^{p-1}h)$ is positive, compact and self-adjoint in $H^1$ with the stated inner product. Compactness follows from the compact embedding into $L^2$ and elliptic inversion. Let $P$ be the span of all nonpositive spectral directions of the Hessian $1-p\lambda L_\Psi^{-1}(u_*^{p-1}\cdot)$. It is finite-dimensional, has smooth basis functions, and its orthogonal complement has a positive spectral gap.

The map $u\mapsto u-\lambda L_\Psi^{-1}(u^p)$ is real analytic near the positive function $u_*$ in $C^{2,\beta}$, $0<\beta<1$. On $P^\perp$ its derivative is invertible: the Hilbert spectral gap gives an inverse in $H^1$, and elliptic Schauder estimates give its boundedness in $C^{2,\beta}$. The analytic implicit function theorem proves \eqref{reduced-equation}. Since $P$ is a spectral subspace, differentiation at $z=0$ gives $D\zeta(0)=0$. Shrinking the parameter neighborhood preserves positivity of $u_z$ and $L_\Psi u_z$. The equation and its finite-dimensional residual give uniform bounds in every fixed smooth norm.

The Hilbert gradient of $I_\lambda$ at $u_z$ belongs to $P$. Since $\partial_{z_\ell}u_z=\phi_\ell+\partial_{z_\ell}\zeta$ with the second term in $P^\perp$, its $P$ coordinates are $\partial_{z_\ell}I_\lambda(u_z)$. Therefore
\begin{equation}\label{background-residual}
L_\Psi u_z-\lambda u_z^p
=\sum_\ell\partial_{z_\ell}I_\lambda(u_z)\,L_\Psi\phi_\ell.
\end{equation}
The classical finite-dimensional \L{}ojasiewicz gradient inequality \cite[(0.2)]{CM15} applies to the analytic function $z\mapsto I_\lambda(u_z)$ and proves \eqref{finite-dimensional-gradient}. Its exponent can be taken in $(1,2]$ on the right by weakening the inequality if necessary.

Differentiate \eqref{reduced-equation}. The Hilbert Hessian at $u_z$ sends $\partial_{z_\ell}u_z$ into $P$. Pairing this assertion with $w\perp P$ gives \eqref{tangent-cancellation}. This identity uses the tangent vector of the reduced family, not merely the fixed basis vector $\phi_\ell$.

Any stationary solution $\widetilde u$ sufficiently close to $u_*$
in $C^{2,\beta}$ belongs to this family. In fact,
$\widetilde u-u_*$ has a $P$-component $\sum z_\ell\phi_\ell$,
and its $P^\perp$-component solves \eqref{reduced-equation}.
Uniqueness in the implicit function theorem gives
$\widetilde u=u_z$. Shrink the neighborhood so the Hessian is
uniformly positive on $P^\perp$ there. A decomposition with
limit $\widetilde u=u_{z_*}$ then has $z_j\to z_*$ using the same
projection and the same family. The modulation and transverse
estimates below hold in this chart because the $P$-projection of
each $\partial_{z_\ell}u_z$ is $\phi_\ell$ and the background
Hessian at $\widetilde u$ is positive on $P^\perp$.
\end{proof}

\begin{lemma}\label{modulation}
For a sequence in Proposition \ref{decomposition}, replace its local bubbles by the family in Proposition \ref{source-convergence}. After small parameter adjustments one can write
\begin{equation}\label{modulation-decomposition}
u_j=u_{z_j}+\sum_{i=1}^k b_{i,j}B_{i,j}+w_j,
\end{equation}
where $z_j\to0$ in a chart based at $u_*$, $b_{i,j}\to(Y/\lambda)^{(n-2)/4}$, and $\|w_j\|_{H^1}\to0$. If $u_*=0$, the background term and $P$ are omitted. In either case, $w_j$ is orthogonal to $P$ and to the amplitude, normalized dilation and normalized translation derivatives of each projected bubble. On this orthogonal complement, for all sufficiently large $j$,
\begin{align}\label{transverse-gap}
\langle w,w\rangle_\Psi
-p\lambda\int_M\left(u_{z_j}+\sum_i b_{i,j}B_{i,j}\right)^{p-1}w^2
\ge c\langle w,w\rangle_\Psi.
\end{align}
\end{lemma}

\begin{proof}
The parameter hypotheses of Proposition \ref{source-convergence}
bound the first two normalized derivatives of each projection
in $H^1$. Each projected mode differs by $o(1)$ in $H^1$
from its localized Euclidean mode. The single-bubble Gram block is consequently the Euclidean Gram block plus $o(1)$, and is invertible. Distinct bubble blocks have inner products tending to zero by \eqref{orthogonal-scales}. Their inner products with fixed smooth background modes also tend to zero. The finite-dimensional implicit function theorem applied to these orthogonality equations gives \eqref{modulation-decomposition}. Bounded second parameter derivatives make the inverse estimates uniform on the normalized parameter neighborhoods.

For completeness, the Euclidean gap can be computed on the round sphere. If $W=U\varphi$ and $\varphi$ is a spherical harmonic of degree $\ell$, then
$$
\frac{c_n\int_{\mathbb R^n}|\nabla W|^2}
{Y\int_{\mathbb R^n}U^{p-1}W^2}
=1+\frac{4\ell(\ell+n-1)}{n(n-2)}.
$$
The removed modes are precisely the degrees zero and one, corresponding to amplitude and to the tangent directions of the family of Sobolev extremals. These are also the directions distinguished in the stability theorem of Bianchi and Egnell \cite{BE91}. The first remaining degree is two, which yields
\begin{equation}\label{euclidean-gap}
c_n\int|\nabla W|^2-pY\int U^{p-1}W^2
\ge\frac4{n+4}c_n\int|\nabla W|^2.
\end{equation}

To transfer this gap, assume \eqref{transverse-gap} fails and normalize $\langle w_j,w_j\rangle_\Psi=1$. Extract its global weak limit and weak limits at every bubble scale. The global limit lies in $P^\perp$, where the background Hessian has a positive gap. Each rescaled limit satisfies the Euclidean orthogonality conditions and obeys \eqref{euclidean-gap}. Their energies have sum at most one. To verify this last assertion even for nested scales, first approximate each rescaled limit by a compactly supported smooth function. Inner products between two resulting rescalings tend to zero when their scales are orthogonal; after summing the finitely many profiles, remove the compact approximation. The same argument includes the global weak limit.

The weighted integrals converge to the sum of the profile integrals.
The $L^{n/2}$ tails of each $U^{p-1}$ are small uniformly outside
large rescaled balls, and local Rellich compactness applies on the
balls. The background term converges by compactness in $L^2$.
By the $H^1$ approximation of the projections, it suffices to
separate the weights of the local Euclidean bubbles. For
$0<\beta<p+1$ the simultaneous conformal change in Lemma
\ref{pair-estimates} takes $\int U_i^{p+1-\beta}U_\ell^\beta$
to $\int_{S^n} b_\mu^\beta$. This tends to zero since $b_\mu\to0$
almost everywhere and $\|b_\mu\|_{p+1}=1$.
Coordinate comparisons give the same conclusion for colliding
centers, and disjoint neighborhoods give it for distinct limiting
centers. Mixed products with the smooth background tend to zero
by subcritical convergence of each bubble. Expanding the
$p-1$ power by the two-term inequality for positive powers,
or bounding its difference by sums of pair minima when $p-1<1$,
shows that the difference from the sum of the weights goes to zero
in $L^{n/2}$. The negative part of the Hessian is therefore at most
$1-c$ times the sum of profile energies, by the two spectral gaps.
This contradicts failure of \eqref{transverse-gap}.
\end{proof}

As in \cite{Bre05}, we establish the energy control along the trajectory.
\begin{proposition}\label{local-gradient}
Let $u_j$ satisfy Proposition \ref{decomposition}. If $u_*=0$, then
\begin{align}\label{pure-gradient}
E_\Psi(u_j)-\lambda\le{}&
C\|L_\Psi u_j-\lambda u_j^p\|_{2n/(n+2)}^2\notag\\
&-c\sum_i\left(\int_Mf_{i,j}L_\Psi^{-1}f_{i,j}
-\frac1{Y}\right)
-c\sum_{i<\ell}\varepsilon_{i\ell,j}.
\end{align}
If $u_*>0$, some $\gamma\in(0,1]$ gives
\begin{align}\label{mixed-gradient}
E_\Psi(u_j)-\lambda\le{}&
C\|L_\Psi u_j-\lambda u_j^p\|_{2n/(n+2)}^{1+\gamma}
-c\sum_i\epsilon_{i,j}^{(n-2)/2}\notag\\
&-c\sum_i\left(\int_Mf_{i,j}L_\Psi^{-1}f_{i,j}
-\frac1{Y}\right)
-c\sum_{i<\ell}\varepsilon_{i\ell,j}.
\end{align}
Empty sums are zero. On a reduced-family neighborhood, the same exponent applies to limiting stationary backgrounds with the same action value as the base point, using the same family with $z_j\to z_*$.
\end{proposition}

\begin{proof}
Use Lemma \ref{modulation}, and temporarily suppress the sequence index. The action comparison will be applied to the function $u-w$ in \eqref{modulation-decomposition}.

Apply \eqref{dual-action} first to the pure source $Y\sum_i b_if_i$. The self-interaction estimate for the chosen family and Lemma \ref{pair-estimates} retain negative contributions of the sizes appearing in \eqref{pure-gradient}. The positive source-overlap errors are smaller, since $n/(n-2)>1$. For each amplitude, direct one-variable maximization gives
\begin{align}\label{amplitude-maximum}
-\frac{Y}2b^2
+\frac p{p+1}\lambda^{-1/p}Y^{(p+1)/p}b^{(p+1)/p}
\le\frac\lambda n\left(\frac{Y}\lambda\right)^{n/2}.
\end{align}
The unique positive maximizer is $b=(Y/\lambda)^{1/(p-1)}$; the expression tends to zero at zero and to minus infinity at infinity. Combining \eqref{amplitude-maximum} with the preceding source estimates gives
\begin{align}\label{pure-action}
I_\lambda(u-w)\le{}&
\frac{k\lambda}n\left(\frac{Y}\lambda\right)^{n/2}
-c\sum_i\left(\int_Mf_iL_\Psi^{-1}f_i-\frac1{Y}\right)\notag\\
&-c\sum_{i<\ell}\varepsilon_{i\ell}.
\end{align}

With a positive background, apply \eqref{dual-action} to $L_\Psi u_z+Y\sum_i b_if_i$. The source remains positive by Lemma \ref{analytic-background}. Its background--bubble source overlap is $O(\sum_i\epsilon_i^{n/2})$. To check the order, $L_\Psi u_z$ has positive upper and lower bounds, and the pointwise source power difference is bounded by $C\min(U_i,U_i^p)$. Splitting at $|x-y_i|\asymp\sqrt{\epsilon_i}$ gives the stated integral. The negative kernel cross-term is at most $-c\sum_i\epsilon_i^{(n-2)/2}$, since $\int_Mu_zf_i\ge c\epsilon_i^{(n-2)/2}$. The background's own Young gap is at most $C\|L_\Psi u_z-\lambda u_z^p\|_2^2$: expand the source power about $\lambda u_z^p$, where the first derivative of this gap vanishes. Equation \eqref{background-residual} then gives
\begin{align}\label{mixed-action}
I_\lambda(u-w)\le{}& I_\lambda(u_z)
+\frac{k\lambda}n\left(\frac{Y}\lambda\right)^{n/2}
+C|\nabla_z I_\lambda(u_z)|^2\notag\\
&-c\sum_i\epsilon_i^{(n-2)/2}
-c\sum_i\left(\int_Mf_iL_\Psi^{-1}f_i-\frac1{Y}\right)
-c\sum_{i<\ell}\varepsilon_{i\ell}.
\end{align}

Pair the true residual $L_\Psi u-\lambda u^p$ with $w$. Orthogonality removes $L_\Psi u_z-\lambda u_z^p$ and the linear bubble sources. Lemma \ref{sum-errors}, the transverse gap, and Taylor expansion yield, in the case of a positive background,
\begin{align}\label{transverse-bound}
\langle w,w\rangle_\Psi\le{}& C\|L_\Psi u-\lambda u^p\|_{2n/(n+2)}^2\notag\\
&+o\left(\sum_i\epsilon_i^{(n-2)/2}
+\sum_i\left(\int_Mf_iL_\Psi^{-1}f_i-\frac1{Y}\right)
+\sum_{i<\ell}\varepsilon_{i\ell}\right).
\end{align}
For a zero background, omit the first sum of scale powers in this formula. The Taylor remainder paired with $w$ is $o(\langle w,w\rangle_\Psi)$: for $p\le2$ it is bounded by $C\|w\|_{p+1}^{p+1}$, and for $p>2$ by $C(\|w\|_{p+1}^3+\|w\|_{p+1}^{p+1})$. These bounds justify absorption. The same orthogonality in the first action variation gives
\begin{align}\label{action-remainder}
I_\lambda(u)\le{}&I_\lambda(u-w)
+C\|L_\Psi u-\lambda u^p\|_{2n/(n+2)}^2\notag\\
&+o\left(\sum_i\epsilon_i^{(n-2)/2}
+\sum_i\left(\int_Mf_iL_\Psi^{-1}f_i-\frac1{Y}\right)
+\sum_{i<\ell}\varepsilon_{i\ell}\right)
\end{align}
in the background case; omit the first sum of scale powers in the pure case.

We next estimate the background gradient without using the square root of the error in \eqref{transverse-bound}. Pair the residual with $\partial_{z_\ell}u_z$. The background contribution is $\partial_{z_\ell}I_\lambda(u_z)$, and the term linear in $w$ for the background operator vanishes by \eqref{tangent-cancellation}. Bubble sources and their nonlinear powers paired with a smooth function contribute $O(\sum_i\epsilon_i^{(n-2)/2})$, because
$$
\int_MB_i+\int_MB_i^p\le C\epsilon_i^{(n-2)/2}.
$$
The change from the background Hessian to the Hessian at $u-w$ costs at most
$$
C\|(u-w)^{p-1}-u_z^{p-1}\|_{2n/(n+2)}\|w\|_{p+1}.
$$
Its squared first factor has the bound \eqref{background-weight-error}. The nonlinear remainder paired with the smooth tangent vector is $O(\langle w,w\rangle_\Psi)$. For $p<2$, this follows from $u-w\ge c>0$, $u\ge0$ and
$$
|u^p-(u-w)^p-p(u-w)^{p-1}w|\le C(u-w)^{p-2}w^2\le Cw^2.
$$
For $p\ge2$, integrate $C((u-w)^{p-2}w^2+|w|^p)$ and apply H\"older's inequality with exponent $p+1$. Consequently, \eqref{transverse-bound} gives
\begin{align}\label{reduced-gradient-bound}
|\nabla_z I_\lambda(u_z)|\le C\biggl(&
\|L_\Psi u-\lambda u^p\|_{2n/(n+2)}
+\sum_i\epsilon_i^{(n-2)/2}\notag\\
&+\sum_i\left(\int_Mf_iL_\Psi^{-1}f_i-\frac1{Y}\right)
+\sum_{i<\ell}\varepsilon_{i\ell}\biggr).
\end{align}
Products of the two square-root error bounds are little-oh of the displayed sums; terms involving the true residual are bounded by a constant times its norm once it is small.

By \eqref{mass-splitting},
$$
I_\lambda(u_*)+\frac{k\lambda}n
\left(\frac{Y}\lambda\right)^{n/2}=\frac\lambda n.
$$
Combine \eqref{mixed-action}, \eqref{action-remainder}, \eqref{finite-dimensional-gradient} and \eqref{reduced-gradient-bound}. The $(1+\gamma)$ and second powers of the positive sums are little-oh of those sums and can be absorbed into their negative terms. The true residual squared is bounded by its $(1+\gamma)$ power for small residual. Finally, mass one gives $E_\Psi(u)-\lambda=2(I_\lambda(u)-\lambda/n)$. This proves \eqref{mixed-gradient}. In the pure case, \eqref{pure-action} and \eqref{action-remainder} give \eqref{pure-gradient} directly. The construction of the reduced family fixes its exponent on its entire parameter neighborhood, as asserted.
\end{proof}

\subsection{Finite length and smooth convergence}\label{finite-length-section}

  We need the following compactness result.
\begin{lemma}\label{stationary-compactness}
Under the hypotheses of Proposition \ref{source-convergence},
for every fixed $\nu>0$ the set of smooth positive solutions
$$
L_\Psi h=\nu h^p,\qquad \int_Mh^{p+1}=1,
$$
is compact in every $C^m(M)$.
\end{lemma}

\begin{proof}
For such solutions $\langle h,h\rangle_\Psi=\nu$, so Proposition \ref{decomposition} applies to any sequence. If a bubble occurs, Proposition \ref{local-gradient}, with $\lambda=\nu$, has zero left side, zero residual and a strictly negative right side. This is impossible. Therefore a subsequence converges strongly in $H^1$ and $L^{p+1}$ to a mass-one positive solution.

Here is the elliptic regularity needed to make this compactness smooth. Strong critical convergence makes the coefficients $\nu h_j^{p-1}$ uniformly integrable in $L^{n/2}$. On small fixed balls their norms are as small as desired. Testing the equation with $\chi^2h_j^{b-1}$ and using the Sobolev inequality gives
\begin{align*}
\|\chi h_j^{b/2}\|_{p+1}^2
\le{}&Cb^2\int_M(|\nabla\chi|^2+1)h_j^b\\
&+Cb^2\|\nu h_j^{p-1}\|_{L^{n/2}(\supp\chi)}
\|\chi h_j^{b/2}\|_{p+1}^2.
\end{align*}
Fix in advance finitely many iteration exponents, starting from $b=p+1$ and ending at some $B>np/2$. Choose the balls so that the last term can be absorbed at every one of these finitely many steps. Shrinking nested cutoffs and a finite covering give a uniform $L^B$ bound. The equation then gives a uniform $W^{2,B/p}$ bound, hence a positive H\"older exponent. A uniformly convergent subsequence has the already identified positive limit, so it has a uniform positive lower bound. Elliptic Schauder estimates now yield bounds in every smooth norm and the required compactness.
\end{proof}

 The following estimate plays important role for the study of the flow asymptotic.
\begin{proposition}
Every trajectory of \eqref{flow} has constants $C>0$, $T>0$ and $\gamma\in(0,1]$ such that
\begin{equation}\label{global-energy-gradient}
0\le r(t)-r_\infty
\le C\|L_\Psi u(t)-r_\infty u(t)^p\|_{2n/(n+2)}^{1+\gamma}
\qquad(t\ge T).
\end{equation}
\end{proposition}

\begin{proof}
For $\lambda=r_\infty$, the possible background masses in
\eqref{mass-splitting} belong to the finite set
$$
 \left\{1-k(Y/\lambda)^{n/2}: k\ge0,\ k(Y/\lambda)^{n/2}\le1\right\}.
$$
For each positive mass $m_*$ in this set, rescaling a background
by $m_*^{-1/(p+1)}$ gives a mass-one stationary solution with
coefficient $\lambda m_*^{2/n}$. Lemma
\ref{stationary-compactness} thus gives compactness of the
backgrounds in $C^{2,\beta}$ on each such mass level.
Cover each level by finitely many neighborhoods from Lemma
\ref{analytic-background}. Every stationary background in a
chosen neighborhood belongs to its reduced family. On that level
the action of a stationary background is $\lambda m_*/n$, so
the reference value of \eqref{finite-dimensional-gradient} is the
same for the base point of the family and the limiting background.
The proof of Proposition \ref{local-gradient} therefore gives
the exponent of the selected family even when the modulation
parameters tend to $z_*\ne0$.

Take $\gamma>0$ to be the minimum of these finitely many exponents
and one. If \eqref{global-energy-gradient} failed, there would be
a sequence $t_j\to\infty$ on which its left side exceeds $j$ times
the residual to the power $1+\gamma$. Propositions
\ref{evolution} and \ref{decomposition} give a subsequence with a
limiting decomposition. In the absence of a positive background,
\eqref{pure-gradient} gives a contradiction. In the presence of
one, select a covering chart of its mass level and apply
\eqref{mixed-gradient} there. Its exponent is at least $\gamma$
and the residual tends to zero, again a contradiction.
\end{proof}

With the understanding above, we can obtain the strong convergence estimate as follows.
\begin{proposition}\label{strong-convergence}
The trajectory converges strongly in $H^1(M)$ and $L^{p+1}(M)$ to one positive smooth solution of \eqref{stationary}. It satisfies
\begin{equation}\label{finite-length}
\int_t^\infty\left(\int_Mu^{p-1}u_t^2\right)^{1/2}\,ds
\le C(r(t)-r_\infty)^{\gamma/(1+\gamma)}
\end{equation}
for all sufficiently large $t$.
\end{proposition}

\begin{proof}
Combining \eqref{global-energy-gradient} and \eqref{residual-dissipation}, and absorbing the term $(r-r_\infty)^{1+\gamma}$ once $r-r_\infty$ is small, gives
\begin{equation}\label{energy-dissipation-gradient}
r-r_\infty\le C\left(\int_Mu^{p-1}u_t^2\right)^{(1+\gamma)/2}.
\end{equation}
If $r-r_\infty$ reaches zero, \eqref{dissipation} makes the trajectory stationary thereafter. Otherwise differentiate its $\gamma/(1+\gamma)$ power and use \eqref{dissipation} and \eqref{energy-dissipation-gradient}:
$$
-\frac d{dt}(r-r_\infty)^{\gamma/(1+\gamma)}
\ge c\left(\int_Mu^{p-1}u_t^2\right)^{1/2}.
$$
Integration proves \eqref{finite-length}. Also,
$$
\left\|\partial_t\left(u^{(p+1)/2}\right)\right\|_2^2
=\frac{(p+1)^2}4\int_Mu^{p-1}u_t^2.
$$
Hence $u^{(p+1)/2}$ has a strong $L^2$ limit of norm one. The elementary inequality
$$
|A^{2/(p+1)}-B^{2/(p+1)}|^{p+1}\le|A-B|^2
\qquad(A,B\ge0)
$$
gives a unique strong $L^{p+1}$ limit for $u(t)$, also of mass one. In every decomposition subsequence this excludes all bubbles. The remaining background is therefore the same nonzero smooth positive function $u_\infty$, and solves \eqref{stationary}. Finally the entire trajectory converges weakly in $H^1$ to $u_\infty$, and
$$
\langle u(t),u(t)\rangle_\Psi=r(t)\longrightarrow r_\infty=\langle u_\infty,u_\infty\rangle_\Psi.
$$
Weak convergence and convergence of the Hilbert norm prove strong $H^1$ convergence.
\end{proof}

Then we can find the convergence result.
\begin{lemma}\label{smoothing}
Assume that a global smooth positive solution of \eqref{flow} satisfies bounded $r(t)$, the lower bound \eqref{curvature-lower}, and strong convergence in $L^{p+1}(M)$ as $t\to\infty$. Then for every $m\ge0$ its tail is bounded in $C^m(M)$, and it has a uniform positive lower bound. The convergence is smooth.
\end{lemma}

\begin{proof}
Strong critical convergence gives uniform nonconcentration in the following sense: if $\rho_j\to0$ and $t_j\to\infty$, then for every fixed $\tau_0>0$,
\begin{equation}\label{nonconcentration}
\sup_{s\in[t_j-\tau_0,t_j]}\sup_{x\in M}
\int_{B_{\rho_j}(x)}u(s)^{p+1}\longrightarrow0.
\end{equation}
This follows by subtracting the strong limit, using the uniform smallness of the tail norm of that difference, and then absolute continuity of the integral of the limiting critical power.

Suppose that $u$ is unbounded. Choose running maximum points $(x_j,t_j)$ with $A_j=u(x_j,t_j)\to\infty$ and $u\le A_j$ on $M\times[0,t_j]$. Smoothness on finite intervals implies $t_j\to\infty$. Set
$$
v_j(y,\tau)=A_j^{-1}u(\exp_{x_j}(A_j^{-2/(n-2)}y),t_j+\tau).
$$
On fixed backward cylinders $0<v_j\le1$ and $v_j(0,0)=1$. Let $g_j$ have the components of $g$ at $A_j^{-2/(n-2)}y$ in the normal coordinates centered at $x_j$. These metrics tend smoothly to the Euclidean metric. Since $p-1=4/(n-2)$, the original time scale is retained, and the equation is
\begin{align}\label{rescaled-parabolic}
\partial_\tau(v_j^p)
={}&\diver_{g_j}\bigl(c_n(v_j^p)^{1/p-1}\nabla(v_j^p)\bigr)\notag\\
&-pA_j^{-4/(n-2)}\Psi(\exp_{x_j}(A_j^{-2/(n-2)}y))v_j
+pr(t_j+\tau)v_j^p.
\end{align}
Applying a smooth convex approximation of the truncation $(v_j^p-2^{-p})_+$ to \eqref{rescaled-parabolic} gives the weak inequality
\begin{align*}
&\partial_\tau(v_j^p-2^{-p})_+\\
&\quad-\diver_{g_j}\bigl(c_n\max\{v_j^p,2^{-p}\}^{1/p-1}
\nabla(v_j^p-2^{-p})_+\bigr)\le C.
\end{align*}
The scalar coefficient is between $c_n$ and $c_n2^{p-1}$. Thus this is an inequality with uniformly bounded ellipticity on the entire cylinder, even though the original diffusion may be singular near zero.

For a radius $R$ independent of $j$, its local mean-value estimate is
\begin{equation}\label{parabolic-mean-value}
(1-2^{-p})\le C
\left(R^{-n-2}\int_{-R^2}^0\int_{B_R}(v_j^p-2^{-p})_+^2\,dv_{g_j}\,d\tau\right)^{1/2}
+CR^2.
\end{equation}
One can obtain this estimate directly as follows. Add the constant $CR^2$ to the truncation, so the right side of its differential inequality is at most $R^{-2}((v_j^p-2^{-p})_++CR^2)$. Testing with $\zeta^2((v_j^p-2^{-p})_++CR^2)^{b-1}$ bounds the time supremum of the squared $L^2$ norm of $\zeta((v_j^p-2^{-p})_++CR^2)^{b/2}$ and its spacetime gradient norm by
$$
Cb^2\iint\bigl(|\nabla\zeta|^2+|\zeta\zeta_\tau|+R^{-2}\zeta^2\bigr)((v_j^p-2^{-p})_++CR^2)^b.
$$
Spatial Sobolev and time H\"older inequalities raise $b$ to $b(1+2/n)$. Iterate from $b=2$ on nested cylinders tending to the half cylinder; the sums of reciprocal exponents and cutoff costs converge. The resulting supremum estimate, evaluated at the terminal time by continuity, is \eqref{parabolic-mean-value}.

Choose $R>0$ so small that the last term in \eqref{parabolic-mean-value} is less than $(1-2^{-p})/2$. The integral of $(v_j^p-2^{-p})_+^2$ on this fixed cylinder is then bounded below by a positive constant. On the other hand,
$$
(v_j^p-2^{-p})_+^2\le(v_j^p)^2\le(v_j^p)^{(p+1)/p}=v_j^{p+1},
$$
because $0\le v_j^p\le1$ and $(p+1)/p<2$. The change of variables gives the exact critical-mass identity
$$
\int_{B_R}v_j(y,\tau)^{p+1}\,dv_{g_j}
=\int_{B_{RA_j^{-2/(n-2)}}(x_j)}u(x,t_j+\tau)^{p+1}\,dv_g.
$$
The spacetime lower bound therefore contradicts \eqref{nonconcentration} on $[t_j-R^2,t_j]$. This proves a uniform upper bound.

The Green-kernel argument in Proposition \ref{evolution}, now using this time-independent upper bound, gives a uniform positive lower bound. We explain why it also gives uniform smooth estimates despite the nonlocal reaction. On every time interval $[t_0,t_0+2]$, make the substitutions
$$
v(t)=\exp\left(-\int_{t_0}^t r(s)\,ds\right)u(t),
\qquad \frac{d\tau}{dt}=\exp\left((1-p)\int_{t_0}^t r(s)\,ds\right).
$$
Both factors have uniform positive upper and lower bounds. The
transformed equation is $v_\tau=c_nv^{1-p}\Delta_gv-\Psi v^{2-p}$,
or $(v^p)_\tau=c_n\diver_g(v^{1-p}\nabla_g(v^p))-p\Psi v$.
The bounds on $v$ give uniform ellipticity and a bounded source,
so \cite[Theorem~6.28]{Lie96} gives an interior H\"older bound for
$v^p$ in local coordinates. Thus the coefficients and source in the
equation for $v$ are H\"older continuous. Interior Schauder estimates
\cite[Theorem~4.9]{Lie96} and iteration give uniform spatial
bounds of all orders on smaller time intervals. These bounds transfer
to $u$, and time derivatives follow from \eqref{flow} and the integral
formula for $r$. Every tail sequence thus has a subsequence converging
in each fixed $C^m$. The strong $L^{p+1}$ limit identifies each
subsequential limit, proving smooth convergence.
\end{proof}

\begin{proof}[Proof of Proposition \ref{source-convergence}]
Proposition \ref{evolution} gives a unique global positive
trajectory and its dissipation identities.
The estimates in Lemmas \ref{pair-estimates} and \ref{sum-errors}
apply to the assumed source family. Proposition \ref{local-gradient}
therefore applies to every concentration subsequence.
Proposition \ref{strong-convergence} gives a unique positive
strong critical limit, and Lemma \ref{smoothing} gives convergence
in every $C^m(M)$.
\end{proof}

\section{Convergence for \texorpdfstring{$\Psi\le R_g$}{potentials at most the scalar curvature}}
\label{all-dimensional-bubbles}

In this section $n\ge4$. To apply Proposition~\ref{source-convergence}
under $\Psi\le R_g$, we compare the projections of the same sources by
the conformal Laplacian and by $L_\Psi$. We adapt Brendle's
correction of the Euclidean bubble \cite{Bre07,BreSurvey} and
the matching to a Green function used by Brendle and Chen for
the boundary Yamabe problem \cite{BC14}. We write $E_{h_y}$
for the geometric Yamabe quotient of a metric $h_y$.

Put $d=\lfloor a\rfloor$ and
$\kappa=((n-2)c_n|S^{n-1}|)^{-1}$.
For the remainder of this section assume $q=R_g-\Psi\ge0$ and $Y_\Psi(M)>0$.
Put $L_g=-c_n\Delta_g+R_g$. Then
$Y(M,[g])\ge Y_\Psi(M)>0$, so both $L_g$ and $L_\Psi$ are coercive.
Their Green kernels are positive and satisfy
$$
c\,d_g(x,z)^{2-n}\le G_g(x,z)\le G_\Psi(x,z)
\le C\,d_g(x,z)^{2-n}\qquad(x\ne z).
$$
Positivity follows by dividing by a positive first eigenfunction
and applying the maximum principle. The fundamental singularity
gives the two-sided bounds near the diagonal; positivity and
compactness give them off it. The middle inequality follows from
$L_\Psi^{-1}-L_g^{-1}=L_\Psi^{-1}qL_g^{-1}$.

All constants in this section may depend on the fixed smooth metric and
on $n$. No assertion is made uniformly over dimensions or over varying
background metrics. We use conformal normal metrics
$h_y=\theta_y^{4/(n-2)}g$ with $\theta_y(y)=1$,
$d\theta_y(y)=0$, and determinant $1+O(r^N)$ in $h_y$-normal
coordinates, where $N>4n+10$.
The conformal normal-coordinate construction
\cite[Theorems~5.1--5.2]{LP87} determines finitely many jets of the
conformal factor by invertible linear maps on homogeneous
polynomials. These maps and their inverses are orthogonally
equivariant. Applying them to the smoothly varying curvature jets,
realizing the resulting jets with a fixed cutoff near the diagonal,
and exponentiating gives a smooth positive family on $M\times M$.
Compactness gives a common coordinate radius and uniform bounds
for its derivatives, including two center derivatives.
The same construction gives $R_{h_y}=O(r^2)$. In normal coordinates at $y$ write
$$
h_y=\exp h,\qquad
h(x)=\sum_{k=2}^{d}H_{k,y}(x)+J_y(x),\qquad
J_y=O''(r^{d+1}).
$$
Here $H_{k,y}$ is homogeneous of degree $k$, is trace-free, and satisfies
$H_{k,y}(x)x=0$. For $n=4,5$ the sum is empty.
The full logarithm satisfies $h(x)x=0$ and $\operatorname{tr}h=O(r^N)$.
Put $t_k(y)=\|H_{k,y}\|_{L^2(S^{n-1})}$. The coordinate construction and
all subsequent operators are orthogonally equivariant, so these definitions
do not require a global orthonormal frame.

Fix a radial smooth function $\chi$ equal to one on $[0,1]$ and zero
on $[2,\infty)$, with $0\le\chi\le1$. Both the polynomial correction
and the final gluing use $\chi(r/\delta)$. We shall choose a single
smooth radius $\delta=\delta(\epsilon)$, independent of the center,
after proving estimates with $\epsilon$ and $\delta$ independent.
For $0<4\epsilon<\delta$ set
$$
\mathcal E_{\epsilon,\delta}(y)
=\epsilon^{n-2}\sum_{k=2}^{d}t_k(y)^2
\int_0^\delta r^{n-1}(\epsilon+r)^{2k+2-2n}\,dr.
$$
Direct integration on $r<\epsilon$ and $\epsilon<r<\delta$ gives
\begin{equation}\label{all-jet-equivalence}
\mathcal E_{\epsilon,\delta}\asymp
\sum_{k<a}t_k^2\epsilon^{2k}
+\boldsymbol1_{\{d=a\}}t_d^2\epsilon^{2a}
\bigl(1+\log(\delta/\epsilon)\bigr).
\end{equation}
The constants in this equivalence are independent of both parameters.

\subsection{The linear correction}

Let $H=\sum_{k=2}^dH_{k,y}$ and $H^\delta=\chi(r/\delta)H$.
Define $w=w_{\epsilon,\delta,y}$ on $\mathbb R^n$ by
\begin{equation}\label{all-linear-correction}
(-c_n\Delta-pYU_\epsilon^{p-1})w
=-U_\epsilon\partial_i\partial_jH^\delta_{ij},
\end{equation}
requiring $w$ to contain only spherical harmonics of degrees at least two
and to have finite Dirichlet energy. Repeated indices are summed.

\begin{lemma}\label{all-linear-lemma}
Equation \eqref{all-linear-correction} has a unique such solution. For
$j=0,1,2$, on $r\le2\delta$ it satisfies
\begin{equation}\label{all-w-pointwise}
|\nabla^jw|\le C\epsilon^a
\sum_{k=2}^dt_k(y)(\epsilon+r)^{k-2a-j}.
\end{equation}
Its spherical mean and first spherical moments vanish on every sphere.
Moreover,
\begin{equation}\label{all-w-energy}
\int_{\mathbb R^n}\left(|\nabla w|^2+
U_\epsilon^{p-1}w^2\right)\le C\mathcal E_{\epsilon,\delta}.
\end{equation}
For sufficiently small fixed $\delta_0$, $U_\epsilon+w\ge U_\epsilon/2$
on $r\le2\delta\le2\delta_0$.

Let $\mathcal Q_{\epsilon,\delta}(H)$ denote minus the coefficient of
$t^2$ in the normalized Yamabe energy of $U_\epsilon+tw$ for the
metric $\exp(tH^\delta)$ on $\mathbb R^n$. Then
\begin{equation}\label{all-quadratic-coercivity}
\mathcal Q_{\epsilon,\delta}(H)\ge c\mathcal E_{\epsilon,\delta}.
\end{equation}
In particular this estimate retains all mixed terms between different
homogeneous degrees.
\end{lemma}

\begin{proof}
If $P$ is trace-free and $P(x)x=0$, then, for every compactly supported
radial function $b$,
$$
\int \partial_i\partial_jP_{ij}\,b(r)=0,\qquad
\int \partial_i\partial_jP_{ij}\,x_\ell b(r)=0.
$$
Indeed the Hessians of $b(r)$ and $x_\ell b(r)$ are sums of multiples
of $\delta_{ij}$, $x_ix_j$, $\delta_{i\ell}x_j$ and
$\delta_{j\ell}x_i$, all of which contract to zero against $P$.
Also
$\partial_i\partial_j(\chi P_{ij})=\chi\partial_i\partial_jP_{ij}$:
the differentiated-cutoff terms vanish by the same contractions and
$x_i\partial_jP_{ij}=-\operatorname{tr}P=0$.
Consequently the right side of \eqref{all-linear-correction} has no
angular components of degree zero or one.

After scaling $x=\epsilon z$, its component of angular degree $\ell\ge2$
is a multiple of
$\epsilon^{k-a-2}r^{k-2}(1+r^2)^{-a}\chi(\epsilon r/\delta)$,
where $\ell\le k-2$. Dividing the operator by $c_n$, the radial operator is
$$
-\partial_r^2-\frac{n-1}{r}\partial_r
+\frac{\ell(\ell+n-2)}{r^2}
-\frac{n(n+2)}{(1+r^2)^2}.
$$
For angular degree one its positive null function is
$r(1+r^2)^{-n/2}$. Dividing a solution by this function turns the
degree-$\ell$ equation into a divergence-form radial equation whose
zeroth-order coefficient is
$[\ell(\ell+n-2)-(n-1)]/r^2>0$.
This gives the maximum principle on finite annuli. The solution regular
at zero is positive and comparable to $r^\ell$ at both ends; the solution
decaying at infinity is positive and comparable to $r^{2-n-\ell}$.
For these assertions, the asymptotic indicial equations at zero and
infinity give the displayed powers. A zero growing coefficient would
give a finite-energy null function, which is impossible by the positive
degree-$\ell$ spherical spectral gap. Their Wronskian therefore gives
a positive Green kernel bounded by
$C r_<^\ell r_>^{2-n-\ell}$.
Inserting the source in this kernel gives, up to a constant,
$$
r^{2-n-\ell}\int_0^r
t^{\ell+n-1}t^{k-2}(1+t^2)^{-a}\,dt
+r^\ell\int_r^{2\delta/\epsilon}
t^{1-\ell}t^{k-2}(1+t^2)^{-a}\,dt.
$$
It is bounded by $C(1+r)^{k-2a}$. The second integral converges at
infinity because $k< n-2+\ell$. At zero the source has the correct
angular order since $k-2\ge\ell$.
Differentiating the integral formula and using its radial equation proves
\eqref{all-w-pointwise}. Outside $2\delta$ the decaying radial solution
gives an additional factor $(2\delta/r)^\ell$ relative to the crude
polynomial bound. Integrating the gradient bound inside $2\delta$ and
this sharper bound outside proves \eqref{all-w-energy}.
The same bounds give $|w|/U_\epsilon\le C\sum t_k(\epsilon+r)^k\le C\delta^2$
on the gluing ball. This proves positivity and the existence assertions.
The exterior decaying solutions give $|w|/U_\epsilon\le C\delta^2$
beyond $2\delta$ as well, so $U_\epsilon+w>0$ on all of $\mathbb R^n$.
After conformal compactification, $w/U_\epsilon$ solves a smooth
homogeneous linear equation near the omitted pole. Its boundedness
there, from the decaying modes, makes the pole removable.

We prove \eqref{all-quadratic-coercivity} without separating the mixed
coefficients. The metric
$\gamma=U_\epsilon^{4/(n-2)}|dx|^2$ is the volume-one round sphere
in stereographic coordinates; let its sectional curvature be $\mathfrak k>0$.
The variation
$$
\gamma_t=(U_\epsilon+tw)^{4/(n-2)}\exp(tH^\delta)
$$
has zero first variation of volume and scalar curvature. The volume
assertion follows from $\int U_\epsilon^pw=0$. For the scalar curvature,
the metric variation of $\Delta U_\epsilon$ is zero because
$H^\delta x=0$ and $\operatorname{tr}H^\delta=0$; conformal covariance
and \eqref{all-linear-correction} give the assertion.

Let $T$ be the transverse trace-free part of $\dot\gamma_0$ on the round
sphere. To justify the decomposition, solve the elliptic equation for
the divergence of the trace-free Lie derivative, on the orthogonal
complement of conformal Killing fields. Its kernel orthogonality follows
by integration by parts. The remainder is $T+f\gamma$ modulo a Lie
derivative. The scalar-curvature variation equation is
$-(n-1)(\Delta_\gamma+n\mathfrak k)f=0$.
Its solutions are the first harmonics and are themselves Lie derivatives
of the round metric. Thus $\dot\gamma_0$ differs from $T$ only by a Lie
derivative. Diffeomorphism invariance and criticality remove all such
terms from the Hessian of the normalized total scalar curvature.
Differentiating the Ricci tensor on a transverse trace-free tensor gives
$\operatorname{Ric}'(T)=\frac12\nabla^*\nabla T+n\mathfrak kT$.
Differentiating the volume normalization then gives
\begin{equation}\label{all-round-hessian}
\mathcal Q_{\epsilon,\delta}(H)
=\frac14\int_{S^n}\left(|\nabla T|^2+2\mathfrak k|T|^2\right)\,dv_\gamma.
\end{equation}

Here is the quantitative step from \eqref{all-round-hessian} to the
polynomial coefficients. Denote by $W'(P)$ the Euclidean linearized
Weyl tensor of a symmetric tensor $P$. On the finite-dimensional space
of polynomials of degrees $2,\ldots,d$ satisfying $Px=0$ and
$\operatorname{tr}P=0$, the map $P\mapsto W'(P)$ is injective.
To verify this, take a homogeneous term of degree $k$. If $W'(P)=0$,
the linearized curvature has the form
$$
R'_{ijkl}=A_{ik}\delta_{jl}+A_{jl}\delta_{ik}
-A_{il}\delta_{jk}-A_{jk}\delta_{il}.
$$
The contracted Bianchi
identity implies $\partial_iA_{jk}=\partial_jA_{ik}$, so
$A=\nabla^2\varphi$ for a homogeneous polynomial $\varphi$ of degree $k$.
In radial gauge the linearized radial Ricci curvature is zero, since
the angular trace of $P$ is zero. Therefore
$$
\Delta\varphi+(n-2)\partial_r^2\varphi=0.
$$
On a degree-$\ell$ spherical component this is
$[(n-1)k^2-\ell(\ell+n-2)]\varphi=0$.
Since $\ell\le k$ and $k\ge2$, the coefficient is positive. Hence
$A=0$. The radial sectional curvature is
$-k(k+1)r^{k-2}P|_{S^{n-1}}/2$, so $P=0$.
Different homogeneous degrees are separated by their polynomial degrees.

This injectivity gives, on the fixed annulus $1<|x|<2$,
$$
\sum_{k=2}^d\|P_k\|_{L^2(S^{n-1})}^2
\le C\|W'(P)\|_{H^{-1}(1<|x|<2)}^2,
\qquad P=\sum_{k=2}^dP_k.
$$
Here the negative norm is dual to $H^1_0$ on the annulus, component
by component. If the inequality failed, a sequence with coefficient
norm one would have a convergent subsequence in this finite-dimensional
space. Its limit would satisfy $W'(P)=0$ on the annulus, hence
everywhere by polynomiality, contradicting injectivity. Thus the
constant controls the entire polynomial, including all mixed degrees,
and depends only on $n$ and $d$. In particular it is independent of
the center and of the two small parameters.
Write the stereographic components of $T$ as
$T=U_\epsilon^{4/(n-2)}\widetilde T$.
Linearized conformal and diffeomorphism invariance of Weyl curvature
give $W'(H^\delta)=W'(\widetilde T)$.
Under $x=\epsilon z$, the round metric pulls back to
$\gamma_1=U_1^{4/(n-2)}|dz|^2$, and the polynomial on the uncut region
becomes $\sum_k\epsilon^kH_{k,y}(z)$. Use the same letters $T$
and $\widetilde T$ for the pulled-back tensor and its components.
On $R<|z|<2R$ with
$R\ge1$ and $2R<\delta/\epsilon$, the further dilation $z=R\xi$
changes its polynomial coefficients to $\epsilon^kR^kH_{k,y}$.
Linearized Weyl curvature is a second-order operator with constant
coefficients, so integration by parts once gives
$$
\|W'(\widetilde T(R\,\cdot))\|_{H^{-1}(1<|\xi|<2)}
\le C\|\widetilde T(R\,\cdot)\|_{H^1(1<|\xi|<2)}.
$$
Applying the coefficient estimate to the full dilated polynomial and
changing variables therefore gives
$$
\sum_k t_k^2\epsilon^{2k}R^{2k}
\le C\left(
R^{2-n}\int_{R<|z|<2R}|\partial\widetilde T|^2\,dz
+R^{-n}\int_{R<|z|<2R}|\widetilde T|^2\,dz\right).
$$
Multiply by $R^{2-n}=R^{-2a}$. Since
$U_1^2\asymp R^{4-2n}$ on this annulus, we obtain
$$
\sum_k t_k^2\epsilon^{2k}R^{2k-2a}
\le C\int_{R<|z|<2R}U_1^2
\left(|\partial\widetilde T|^2+|z|^{-2}|\widetilde T|^2\right)\,dz.
$$

To sum the right sides, use the conformal metric
$\gamma_1=U_1^{4/(n-2)}|dz|^2$. Its connection differs from the
Euclidean connection by $O(|z|^{-1})$ for $|z|\ge1$, and
$|T|_{\gamma_1}=|\widetilde T|$. Consequently the gradient integrals
are bounded by $C\|\nabla T\|_{L^2(\gamma_1)}^2$ plus the
displayed zeroth-order integrals. The latter satisfy
$$
\int_{|z|>1}U_1^2|z|^{-2}|\widetilde T|^2\,dz
\le C\int_{S^n}
\operatorname{dist}_{\gamma_1}(x,\mathrm{pole})^{-2}
|T|_{\gamma_1}^2\,dv_{\gamma_1}
\le C\|T\|_{H^1(\gamma_1)}^2.
$$
The first inequality uses
$dv_{\gamma_1}=U_1^{2n/(n-2)}\,dz$ and
$\operatorname{dist}_{\gamma_1}(x,\mathrm{pole})\asymp |z|^{-1}$
near the omitted pole. For the second, apply the scalar local Hardy
inequality to $|T|_{\gamma_1}$ and then Kato's inequality. The Hardy
inequality follows by integrating
$\operatorname{div}(x/|x|^2)=(n-2)/|x|^2$ in a fixed coordinate ball
and applying Cauchy's inequality; the derivative of its fixed cutoff
contributes only an $L^2$ term.

Take the disjoint annuli with $R=1,2,4,\ldots$ and
$2R<\delta/\epsilon$. For $k<a$ the geometric sum of
$R^{2k-2a}$ is bounded above and below by positive constants.
For $k=a$ the number of annuli is comparable to
$1+\log(\delta/\epsilon)$. The first annulus also handles
$4<\delta/\epsilon<8$. Thus the sum of the left sides controls
$\mathcal E_{\epsilon,\delta}$ by \eqref{all-jet-equivalence},
whereas the sum of the right sides is controlled by
$\|T\|_{H^1(\gamma_1)}^2$. The dilation is an isometry between
$\gamma$ and $\gamma_1$. Equation \eqref{all-round-hessian}
therefore proves \eqref{all-quadratic-coercivity}.
\end{proof}

\subsection{Green matching and the mass term}\label{all-mass-section}

Let $G_y=\kappa^{-1}G_{h_y}(\cdot,y)$, so that $r^{n-2}G_y\to1$.
It is positive. The local Green estimate needed for matching is
\begin{equation}\label{all-green-bound}
|\nabla^j(G_y-r^{-2a})|
\le C\sum_{k=2}^dt_k(y)r^{k-2a-j}+Cr^{d+1-2a-j},
\qquad j=0,1,2.
\end{equation}
For completeness, applying $L_{h_y}$ to a cut-off $r^{-2a}$ leaves,
near zero, a source bounded by
$C\sum t_kr^{k-n}+Cr^{d+1-n}$; the radial Laplacian error from the
determinant is of higher order. Convolution with the bound
$C|x-z|^{2-n}$ for the Green kernel gives
$Cr^{k+2-n}$ for each exponent $0<k<n-2$.
This follows by splitting $|z|<r/2$, $|x-z|<r/2$, and their complement
and integrating the corresponding radial powers. For $n\ge5$, all these exponents satisfy
$k\le d+1<n-2$. The nonsingular part is bounded and is included
in the last term. Interior estimates after rescaling each annulus
give the two derivative bounds.

In dimension four we use the conformal normal-coordinate expansion
\cite[Lemma 6.4]{LP87}, which gives
$G_y=r^{-2}+A_y+O''(r)$, with uniform annular bounds.
These bounds follow also from the smooth dependence of the finite
conformal jets and the local Green parametrix. Since $d+1-2a=0$
in dimension four, they give \eqref{all-green-bound}.

Define
\begin{align*}
\mathcal I_y(\delta)
&=c_n\int_{|x|=\delta}
\left(r^{-2a}\partial_rG_y-G_y\partial_rr^{-2a}\right)\,d\sigma,\\
\mathcal Z&=\{y:H_{2,y}=\cdots=H_{d,y}=0\}.
\end{align*}
The measure in this definition is Euclidean. Replacing it by the normal
coordinate measure costs $O(\delta^{N+2-n})$ in the estimates below.

\begin{lemma}\label{all-mass-flux}
On $\mathcal Z$ the limit $\mathcal I_y(0)$ exists uniformly and is
continuous in $y$, and
$$
\mathcal I_y(\delta)=\mathcal I_y(0)+O(\delta^{\tau}),
\qquad \tau=2d+4-n>0.
$$
In a nonround positive Yamabe conformal class it is strictly positive.
\end{lemma}

\begin{proof}
At $y\in\mathcal Z$, $h=O''(r^{d+1})$ and
$R_{h_y}=\partial_i\partial_jh_{ij}+O(r^{2d})+O(r^{N-2})$.
The spherical average of the linear term is $O(r^{N-2})$, by the
radial-gauge integration by parts used above. Also
$G_y-r^{-2a}=O''(r^{d+1-2a})$.
Averaging $L_{h_y}G_y=0$, the angular Laplacian integrates to zero
up to the higher determinant error. Thus, if a bar denotes spherical
average,
$$
\frac d{dr}\mathcal I_y(r)
=|S^{n-1}|r\,\overline{R_{h_y}G_y}
+O(r^{N-n+1})=O(r^{2d+1-2a}).
$$
The exponent is $\tau-1>-1$. Integration proves uniform convergence;
continuity follows because $\mathcal I_y(r)$ is continuous for each
fixed positive $r$.

For fixed $y\in\mathcal Z$, the Green metric
$g_*=G_y^{4/(n-2)}h_y$ is complete, smooth, scalar flat and one-ended.
In the inverted coordinates with $\rho=r^{-1}$, the metric error and
its derivatives of orders $j=0,1,2$ are $O(\rho^{-d-1-j})$.
The normal metric contribution to the ADM flux tends to zero.
After inversion its error $\eta$ satisfies
$\eta x=0$ and
$\operatorname{tr}\eta=O(\rho^{-2d-2})+O(\rho^{-N})$.
The products of the metric and conformal errors have flux
$O(\rho^{n-2-2d-2})=O(\rho^{-\tau})$.
The remaining conformal flux, using the usual ADM normalization, is
$$
\mathfrak m_{\rm ADM}
=\frac{\mathcal I_y(0)}{a c_n|S^{n-1}|}.
$$
For example, if $G_y=r^{-2a}+A+o(1)$ in spherical mean, this gives
$\mathfrak m_{\rm ADM}=2A$, since
$\mathcal I_y(0)=2ac_n|S^{n-1}|A$.

For orientable $M$ and $4\le n\le7$, take the time-symmetric data
with metric $g_*$ in \cite[Theorem~1 and Section~2]{EHLS16}.
The derivative bounds give the weighted Sobolev decay there for
integrability exponents above $n$ and rates in
$(a,\min\{d+1,n-2\})$. Both constraint densities and the ADM
momentum vanish, so the mass is nonnegative.
The derivative bounds also give the weighted H\"older decay for
the equality argument in \cite[Lemma~10.7]{LP87}. Its scalar-flat
conformal corrections retain admissible decay at a slightly smaller
rate, so the same mass inequality yields Euclidean rigidity.

For $n\ge8$ we use \cite[Theorem~1.1 and Definition~2.2]{KWW26}.
The metric $g_*$ is complete, has empty singular set, and has
order $d+1>(n-2)/2$ with the stated derivative bounds. Its scalar
curvature is zero and hence integrable. The theorem gives both
nonnegative mass and Euclidean rigidity.

A Euclidean Green metric conformally compactifies to the round
sphere. For nonorientable $M$, lift $h_y$ and a positive first
eigenfunction of $L_{h_y}$ to the orientable double cover.
The lifted eigenfunction gives positive Yamabe type on the cover.
The conformal normal jets lift at both poles, and the lift of
$G_y$ is the sum of their normalized single-pole Green functions.
The other summand adds a positive contribution to the nonnegative
flux at either pole. This proves strict positivity in a nonround class.
\end{proof}

Set
$$
V_{\epsilon,\delta,y}
=\chi(r/\delta)(U_\epsilon+w_{\epsilon,\delta,y})
+(1-\chi(r/\delta))K\epsilon^aG_y.
$$
It is defined by the Green term outside the coordinate ball. It is
positive and bounded above and below by constant multiples of the
bubble profile, with constants independent of the two parameters.

\begin{lemma}\label{all-two-parameter}
There are constants $c,C,\delta_0>0$, independent of the center,
such that for $0<4\epsilon<\delta<\delta_0$,
\begin{align}\label{all-energy-cost}
\int_M V L_{h_y}V-Y\left(\int_MV^{p+1}\right)^{2/(p+1)}
\le{}&-c\mathcal E_{\epsilon,\delta}\notag\\
&\quad
-K^2\epsilon^{2a}\mathcal I_y(\delta)\notag\\
&+C\epsilon^{2a}\sum_{k=2}^dt_k\delta^{k-2a}
+C\epsilon^{2a}\delta^\tau+C\epsilon^n\delta^{-n}.
\end{align}
For $v=V/\|V\|_{p+1}$,
\begin{equation}\label{all-residual-two-limit}
\limsup_{\epsilon\downarrow0}\sup_{y\in M}
\frac{\|L_{h_y}v-Yv^p\|_s^2}
{\mathcal E_{\epsilon,\delta}(y)+\epsilon^{2a}}
\le C\delta^\tau
\end{equation}
for each fixed $\delta<\delta_0$. The scale limits in this lemma are
also uniform when $\delta$ lies in a compact subinterval of
$(0,\delta_0)$. All integrals and norms in this
lemma use $h_y$.
\end{lemma}

\begin{proof}
We give the energy and equation estimates separately, retaining the
same correction and cutoffs. First extend $\chi(r/\delta)h$ by zero
on Euclidean space. Up to errors of order $r^N$ it remains trace-free
and radial. For its low-degree part $H^\delta$, Taylor expansion of
the normalized energy at $(|dx|^2,U_\epsilon)$ has zero linear term
and quadratic term $-\mathcal Q_{\epsilon,\delta}(H)$.
The following estimates justify the expansion and its use with the
remaining higher-degree terms.

The Christoffel formula for $\exp h$ gives
$$
R_{\exp h}=\partial_i\partial_jh_{ij}
+R^{(2)}(h)
+O(|h|^2|\partial^2h|+|h||\partial h|^2)
$$
when the trace is zero; the determinant error is $O(r^{N-2})$.
Here $R^{(2)}$ is exactly the quadratic Taylor coefficient. Integrating
the term with $\partial^2h$ by parts against $U_\epsilon^2$ bounds
the cubic error by a constant times
$$
\|h\|_\infty\int U_\epsilon^2
\left(|\partial h|^2+(\epsilon+r)^{-2}|h|^2\right).
$$
The differentiated weights satisfy
$|\partial U_\epsilon|\le C U_\epsilon/(\epsilon+r)$.
The cutoff derivatives on the annulus obey the same weighted bound.
The scalar cross-remainders obey the same estimate. For example,
integrating $R^{(1)}(h)w^2$ once by parts gives the bound
$$
C\|w/U_\epsilon\|_\infty
\left(\int U_\epsilon^2|\partial h|^2\right)^{1/2}
\|\nabla w\|_2.
$$
In $R^{(2)}(h)U_\epsilon w$, integrate each $h\partial^2h$ term
once by parts. The resulting terms have either the factor
$\|w/U_\epsilon\|_\infty$ times the weighted quadratic metric norm,
or $\|h\|_\infty$ times its square root times $\|\nabla w\|_2$.
They are bounded by $C\delta^2\mathcal E_{\epsilon,\delta}$.
The radial property makes the metric correction to
$|dU_\epsilon|^2+2\langle dU_\epsilon,dw\rangle$ vanish;
the remaining Dirichlet correction is bounded by
$C\|h\|_\infty\int|dw|^2$.
Since $|w|/U_\epsilon\le C\delta^2$, the cubic error in the mass
normalization is bounded by
$C\delta^2\int U_\epsilon^{p-1}w^2$.
These estimates and \eqref{all-w-energy} bound the nonlinear remainder
for $H^\delta$ by $C\delta^2\mathcal E_{\epsilon,\delta}$.
The same bound applies to the passage from the normalized energy to
the difference of the numerator and $Y$ times the critical norm squared:
the mass change is $O(\mathcal E_{\epsilon,\delta})$ and
$\mathcal E_{\epsilon,\delta}\le C\delta^4$.

The higher part $J=O''(r^{d+1})$ has zero linear pairing with every
radial weight, apart from the higher determinant error. Its weighted
quadratic norm is at most
$$
C\epsilon^{2a}\int_0^{2\delta}
r^{n-1}(\epsilon+r)^{2d+4-2n}\,dr
\le C\epsilon^{2a}\delta^\tau.
$$
Each quadratic cross-term with $H^\delta$ or $w$ is bounded by
$\eta\mathcal E_{\epsilon,\delta}+C_\eta\epsilon^{2a}\delta^\tau$,
using Cauchy's inequality before estimating the coefficients.
In these pairings the additional norm
$\int w^2/(\epsilon+r)^2$ is bounded by the Euclidean Hardy inequality
and \eqref{all-w-energy}.
Integration by parts in $R^{(2)}$ uses a compactly supported extension,
so no boundary term is lost. The determinant errors are bounded by
$C\epsilon^{2a}\delta^{N+2-n}$ and are smaller than the displayed
higher-jet cost. Thus the full Euclidean extension has energy deficit
at least $c\mathcal E_{\epsilon,\delta}-C\epsilon^{2a}\delta^\tau$,
after fixing $\eta$ and then $\delta_0$.

We now replace its exterior by the Green tail. In the following
identities write $\chi=\chi(r/\delta)$. Off the pole, the product
rule and $L_{h_y}G_y=0$ give the exact decomposition
\begin{align*}
L_{h_y}V
={}&\chi L_{h_y}(U_\epsilon+w)
-c_n(\Delta_{h_y}\chi)(U_\epsilon+w-K\epsilon^aG_y)\\
&-2c_n\langle d\chi,d(U_\epsilon+w-K\epsilon^aG_y)\rangle_{h_y}.
\end{align*}
Both differentiated-cutoff terms are supported on
$\delta<r<2\delta$, where
$|d\chi|\le C\delta^{-1}$ and
$|\Delta_{h_y}\chi|\le C\delta^{-2}$.
On that annulus \eqref{all-green-bound} and
\eqref{all-w-pointwise} give
\begin{align*}
|V-K\epsilon^aG_y|+\delta^2|L_{h_y}V|
\le C\epsilon^a\left(
\sum_kt_k\delta^{k-2a}+\delta^{d+1-2a}
+\epsilon^2\delta^{-n}\right).
\end{align*}
For $r\ge2\delta$ one has $L_{h_y}V=0$.
The exact Green identity on $M\setminus B_\delta$ is
\begin{align*}
\int_{M\setminus B_\delta}\left(c_n|dV|_{h_y}^2+R_{h_y}V^2\right)
+c_n\int_{\partial B_\delta}V\partial_rV
={}&\int_{M\setminus B_\delta}
(L_{h_y}V)(V-K\epsilon^aG_y)\\
&-c_nK\epsilon^a\int_{\partial B_\delta}
(V\partial_rG_y-G_y\partial_rV).
\end{align*}
All integrals in this identity use the actual metric and induced boundary
measure. This is obtained by integrating the two products with $G_y$, for which
$L_{h_y}G_y=0$ off the pole; the normal on the exterior domain is
$-\partial_r$.
The bulk integrand is supported on the matching annulus.
Multiplying the preceding pointwise bounds and using its volume
$O(\delta^n)$ bounds its integral by
$$
C\epsilon^{2a}\sum_kt_k^2\delta^{2k-2a}
+C\epsilon^{2a}\delta^\tau+C\epsilon^n\delta^{-n}.
$$
For the higher-jet term the power is
$2(d+1-2a)+n-2=\tau$; the pure bubble-tail error gives
$\epsilon^{n+2}\delta^{-n-2}\le C\epsilon^n\delta^{-n}$.
Cross-terms are bounded by the corresponding squares.
In the Wronskian replace $V|_{r=\delta}=U_\epsilon+w$ by
$K\epsilon^ar^{-2a}$. The $w$ contribution is bounded by
$C\epsilon^{2a}\sum_kt_k\delta^{k-2a}$.
The difference $U_\epsilon-K\epsilon^ar^{-2a}$ has two derivative
bounds $C\epsilon^{a+2}r^{-2a-2-j}$ and costs
$C\epsilon^n\delta^{-n}$. The main Wronskian is exactly
$-K^2\epsilon^{2a}\mathcal I_y(\delta)$. Its sign follows from
the exterior normal $-\partial_r$ in the preceding identity.
With Euclidean surface measure, the leading singular part of the
boundary pairing with $w$ has zero integral, since both $w$ and
its radial derivative have zero spherical mean. The change to
the metric surface measure contributes only the determinant
error already estimated. The stated bound controls the pairing
with $G_y-r^{-2a}$.

To see explicitly which exterior energy is removed, put
$\bar h=\exp(\chi h)$ and $W=U_\epsilon+w$ on $\mathbb R^n$.
The Euclidean equation for $U_\epsilon$ gives
$$
c_n\int_{r>\delta}|dU_\epsilon|^2
+c_n\int_{r=\delta}U_\epsilon\partial_rU_\epsilon
=Y\int_{r>\delta}U_\epsilon^{p+1}
=O(\epsilon^n\delta^{-n}).
$$
In the corresponding expression with $\bar h,W$, every term linear
in $w$ vanishes by its spherical mean zero, including the boundary
pairing. The linear scalar-curvature term has zero spherical mean
by radial gauge and trace zero, up to the stated determinant error.
The remaining terms involve $|dw|^2$, $R^{(2)}(\chi h)U_\epsilon^2$
or $(\partial^2(\chi h))U_\epsilon w$. The pointwise bounds and one
integration by parts bound their sum, including the quadratic
boundary term, by
$$
C\epsilon^{2a}\sum_kt_k^2\delta^{2k-2a}
+C\epsilon^{2a}\delta^\tau+C\epsilon^n\delta^{-n}.
$$
Thus subtraction of this exterior expression from the preceding Green
identity leaves precisely its negative flux and these costs.
Both exterior critical masses are $O(\epsilon^n\delta^{-n})$.
Since the coefficients $t_k$ are bounded and $\delta$ is small,
$t_k^2\delta^{2k-2a}\le C t_k\delta^{k-2a}$.
Combining these estimates proves \eqref{all-energy-cost}.

For the equation estimate, introduce only in this calculation
$$
N_m=\|U_\epsilon(\epsilon+r)^{m-2}\|_{L^s(r<2\delta)}.
$$
Integration of
$\epsilon^{as}r^{n-1}(\epsilon+r)^{s(m-n)}$ gives
$$
N_m\le C\begin{cases}
\epsilon^m,&m<a,\\
\epsilon^a(1+\log(\delta/\epsilon))^{1/s},&m=a,\\
\epsilon^a\delta^{m-a},&m>a.
\end{cases}
$$
On $r<\delta$, where $H^\delta=H$, the Euclidean equations give
the exact identity
\begin{align*}
L_{h_y}(U_\epsilon+w)-Y(U_\epsilon+w)^p
={}&-c_n(\Delta_{h_y}-\Delta)(U_\epsilon+w)\\
&+(R_{h_y}-\partial_i\partial_jH_{ij})U_\epsilon+R_{h_y}w\\
&-Y\bigl((U_\epsilon+w)^p-U_\epsilon^p
-pU_\epsilon^{p-1}w\bigr).
\end{align*}
The radial gauge and determinant estimate control the Laplacian
error on $U_\epsilon$; the linear correction has cancelled the
low-degree scalar-curvature term. The remaining terms are bounded
by a constant times
$$
\sum_{i,j=2}^dt_it_jU_\epsilon(\epsilon+r)^{i+j-2}
+U_\epsilon(\epsilon+r)^{d-1}.
$$
For the metric error this follows from
$|h||\partial^2w|+|\partial h||\partial w|$ and the scalar-curvature
expansion above. For the nonlinear term use
$|w|\le U_\epsilon/2$ and
$U_\epsilon^{p-1}\le C\epsilon^2(\epsilon+r)^{-4}$.
The determinant term has a higher power and is included in the second
summand. Its norm is therefore at most
$C\sum t_it_jN_{i+j}+C\epsilon^a\delta^{d+1-a}$.

On the matching annulus use the product-rule identity above.
The nonlinear term there is bounded by $C U_\epsilon^p$.
Multiplying the pointwise bounds by
$\delta^{n/s}=\delta^{(n+2)/2}$ gives norm at most
$$
C\epsilon^a\left(\sum_kt_k\delta^{k-a}
+\delta^{d+1-a}+\epsilon^2\delta^{-a-2}\right).
$$
In the exterior $L_{h_y}V=0$, and the remaining nonlinear term has
norm at most $C\epsilon^{a+2}\delta^{-a-2}$, by integration of the
Green upper bound.
For fixed $\delta$, the square of each low-degree term is
$o(\mathcal E_{\epsilon,\delta})$, uniformly over the bounded
coefficient set. For example, if $i\le j$ and $i+j<a$, dividing
$(t_it_jN_{i+j})^2$ by $t_i^2\epsilon^{2i}$ leaves
$Ct_j^2\epsilon^{2j}$. If $i+j=a$, the additional factor is
$(1+\log(\delta/\epsilon))^{2/s}$ and the limit is still zero.
If $i+j>a$ and $i<a$, the remaining bound is
$Ct_j^2\epsilon^{2(a-i)}\delta^{2(i+j-a)}$.
If $i=a$, then $i=j=a$ and division by the critical term leaves
$Ct_a^2\delta^{2a}/(1+\log(\delta/\epsilon))$.
The annular terms use the factors $(\epsilon/\delta)^{2(a-k)}$
for $k<a$, and $(1+\log(\delta/\epsilon))^{-1}$ for $k=a$.
The higher-degree terms have square at most
$C\epsilon^{2a}\delta^{2(d+1-a)}=C\epsilon^{2a}\delta^\tau$.

Finally spherical mean zero of $w$, its quadratic mass bound, and the
exterior estimate give
$$
\left|\int_MV^{p+1}-1\right|
\le C\mathcal E_{\epsilon,\delta}+C\epsilon^n\delta^{-n}.
$$
The square of this normalization error is
$o_\delta(\mathcal E_{\epsilon,\delta}+\epsilon^{2a})$.
Indeed $\mathcal E_{\epsilon,\delta}\to0$ uniformly in the
bounded coefficient set for fixed $\delta$, and
$\epsilon^{2n}\delta^{-2n}/\epsilon^{2a}\to0$.
The exact normalization identity is
$$
L_{h_y}v-Yv^p
=\frac{L_{h_y}V-YV^p}{\|V\|_{p+1}}
+Y\bigl(\|V\|_{p+1}^{-1}-\|V\|_{p+1}^{-p}\bigr)V^p.
$$
Since $\|V^p\|_s=\|V\|_{p+1}^p\to1$, this proves
\eqref{all-residual-two-limit}.
All constants used in these regional estimates are uniform for
$y\in M$ and for $\delta$ in a compact subinterval of
$(0,\delta_0)$. The scale limit is taken with $\delta$ fixed;
the next proposition makes the joint choice of $\delta$ and
$\epsilon$ needed for a single smooth family.
\end{proof}

\subsection{The family of projected sources}

\begin{proposition}\label{all-positive-source}
Suppose $Y(M,[g])>0$ and the conformal class is not round. There is
one smooth function $\delta(\epsilon)\downarrow0$, independent of $y$,
such that the normalized positive functions
$$
v^g_{\epsilon,y}=\theta_y
\frac{V_{\epsilon,\delta(\epsilon),y}}
{\|V_{\epsilon,\delta(\epsilon),y}\|_{L^{p+1}(h_y)}},\qquad
f_{\epsilon,y}=(v^g_{\epsilon,y})^p,
\qquad B^g_{\epsilon,y}=YL_g^{-1}f_{\epsilon,y}
$$
satisfy, uniformly for all centers,
\begin{gather}\notag
\|f_{\epsilon,y}\|_s=1,\qquad
\int_Mf_{\epsilon,y}\ge c\epsilon^a,\qquad
0<B^g_{\epsilon,y}(x)\le
C\left(\frac{\epsilon}{\epsilon^2+d_g(x,y)^2}\right)^a,\\
\label{all-dual-properties}
D_g(\epsilon,y):=\int_Mf_{\epsilon,y}L_g^{-1}f_{\epsilon,y}-Y^{-1}
\ge c\left(\mathcal E_{\epsilon,\delta(\epsilon)}(y)+\epsilon^{n-2}\right),
\qquad D_g\longrightarrow0,\\
\label{all-small-residual}
\|L_gB^g_{\epsilon,y}-Y(B^g_{\epsilon,y})^p\|_s^2=o(D_g(\epsilon,y)).
\end{gather}
The source and its first two normalized scale and center derivatives
are bounded in $L^s$. The corresponding derivatives of $B^g$ are
bounded in $H^1$ and differ by $o(1)$ from the respective local
Euclidean modes.
On a round conformal class the exact spherical family specified
at the end of the proof has the same profile and derivative
conclusions, and $D_g=0$ with zero geometric residual.
\end{proposition}

\begin{proof}
First hold $\delta$ fixed. Each linear annular cost in
\eqref{all-energy-cost} is absorbed by a small fraction of
$\mathcal E_{\epsilon,\delta}$ plus $o_\delta(\epsilon^{2a})$.
For $k<a$, Young's inequality leaves
$C_\eta\epsilon^{4a-2k}\delta^{2k-4a}=o_\delta(\epsilon^{2a})$.
For $k=a$ it leaves
$$
C_\eta\epsilon^{2a}\delta^{-2a}/(1+\log(\delta/\epsilon))
=o_\delta(\epsilon^{2a}).
$$
If $\mathcal Z\ne\varnothing$, compactness and Lemma
\ref{all-mass-flux} give a positive lower bound for
$\mathcal I_y(\delta)$ on $\mathcal Z$, simultaneously for all sufficiently
small $\delta$. For each fixed such $\delta$ it persists on a neighborhood
of $\mathcal Z$. Outside that neighborhood $\sum t_k^2$ has a positive
minimum, and \eqref{all-jet-equivalence} gives
$\mathcal E_{\epsilon,\delta}/\epsilon^{2a}\to\infty$ uniformly.
It absorbs the bounded flux of either sign there.
If $\mathcal Z=\varnothing$, the latter argument applies everywhere.
For $n=4,5$ the sum is empty and $\mathcal Z=M$.
Consequently, for all sufficiently small fixed $\delta$, there is an
$\epsilon_\delta>0$ such that
\begin{equation}\label{all-uniform-test-gap}
Y-E_{h_y}(V_{\epsilon,\delta,y})
\ge c_1\mathcal E_{\epsilon,\delta}(y)+c_2\epsilon^{2a}
\end{equation}
for every $y$ and $0<\epsilon<\epsilon_\delta$.
The constants $c_1,c_2>0$ can be chosen independent of small $\delta$.
This argument also works uniformly when $\delta$ ranges over any closed
interval of positive sufficiently small radii.

Here is a joint smooth choice, so no fixed-center limit is being used.
Take $\delta_j=2^{-j}\delta_0$. On the compact interval
$[\delta_{j+1},\delta_j]$, choose $\epsilon_j$ so that
\eqref{all-uniform-test-gap} holds and the residual ratio in
\eqref{all-residual-two-limit} is at most
$C\delta_j^\tau+1/j$ for every center and every $0<\epsilon\le\epsilon_j$.
Decrease $\epsilon_j$ further so that
$\epsilon_{j+1}<\epsilon_j/16$ and $\epsilon_j<\delta_j^{j}$.
Interpolate $\delta$ from $\delta_j$ to $\delta_{j+1}$ as
$\log(1/\epsilon)$ runs from $\log(1/\epsilon_j)$ to
$\log(1/\epsilon_{j+1})$, using a fixed smooth step flat at both ends.
Then
\begin{equation}\label{radius-bounds}
\delta(\epsilon)\to0,\qquad
\epsilon\delta(\epsilon)^{-m}\to0\quad(m\text{ fixed}),\qquad
|(\epsilon\partial_\epsilon)^\ell\delta|\le C_\ell\delta
\quad(\ell=1,2).
\end{equation}
The stated energy gap and the residual limit now hold along this one
radius, uniformly on $M$.

The normalized function $v^g$ has critical norm one, so its source has
$L^s$ norm one exactly. On $r<\epsilon$, positivity and
\eqref{all-w-pointwise} give $\int f\ge c\epsilon^a$.
The profile upper bound for $v^g$ holds on the core, on the annulus,
and on the exterior Green tail separately. The same bound for $B^g$
follows by applying the Green-kernel bound to $f\le CU_\epsilon^p$.
For a point at distance $r>2\epsilon$, the source ball of radius $r/2$
contributes $C\epsilon^ar^{-2a}$; the ball of radius $r/2$ about that
point contributes $C\epsilon^{a+2}r^{-2a-2}$; the remaining dyadic
annuli have a summable radial integral with the same upper bound.
For $r\le2\epsilon$, scaling gives $C\epsilon^{-a}$.

Write $v=v^g$ in the remaining identities. Positivity of $L_g$ and
Cauchy's inequality for its quadratic form give the exact comparison
$$
1=\left(\int f v\right)^2
\le\left(\int fL_g^{-1}f\right)\left(\int vL_gv\right),
\qquad
D_g\ge\frac1{E_g(v)}-\frac1Y.
$$
Thus \eqref{all-uniform-test-gap} proves \eqref{all-dual-properties}.
Furthermore,
$$
B^g-v=-L_g^{-1}(L_gv-Yv^p).
$$
Elliptic inversion from $L^s$ to $H^1$ and the Lipschitz bound for
the power map on bounded subsets of $L^{p+1}$ yield
$$
\|L_gB^g-Y(B^g)^p\|_s
=Y\|v^p-(B^g)^p\|_s\le C\|L_gv-Yv^p\|_s.
$$
This proves \eqref{all-small-residual}. The same identity implies
$B^g-v\to0$ in $H^1$; since $\int v^{p+1}=1$,
$\int f B^g\to1$, which proves $D_g\to0$.

We verify the parameter estimates with the radius already chosen.
For each angular component of the term of degree $k$, the formula
for $w$ has the form
$$
 \epsilon^{k-a}\,
 \mathcal L_\ell^{-1}\!\left[
   F_{k,\ell,y}\,\chi(\epsilon|\cdot|/\delta(\epsilon))
 \right](x/\epsilon),
$$
where $\mathcal L_\ell$ is the fixed rescaled radial operator and
$F_{k,\ell,y}$ is a smooth multiple in $y$ of the uncut radial
source in the proof of Lemma \ref{all-linear-lemma}. A scale
derivative acts on the prefactor, on the argument $x/\epsilon$,
or on the rescaled cutoff. Up to two such cutoff derivatives are
bounded on $\delta/\epsilon<|x/\epsilon|<2\delta/\epsilon$ by
\eqref{radius-bounds}. The radial Green formula therefore gives
\eqref{all-w-pointwise} for these differentiated sources as well,
including the argument derivatives with their factors
$|x|/\epsilon$. In a smooth center frame, $\epsilon\partial_{y^j}$
acts on a coefficient with a factor $\epsilon$, or on $x/\epsilon$
with a bounded coefficient. The same argument applies to ordered
products of two normalized derivatives. In particular their
core bounds are $C U_\epsilon(\epsilon+|x|)^2$ for the correction,
where the coefficient norms and their first two center derivatives
are bounded on the fixed compact manifold.

For the tail one has
$G_{h_y}(x,y)=\theta_y(x)^{-1}G_g(x,y)$ and
$$
 |\partial_y^j G_g(x,y)|\le C d_g(x,y)^{2-n-j}\qquad (0\le j\le2)
$$
by the local Green parametrix and elliptic estimates off the
diagonal. On the matching annulus each normalized center
derivative of the cutoff costs at most $C\epsilon/\delta$,
and each normalized scale derivative is bounded. Hence $V$ and
all its normalized derivatives up to order two are bounded by
the uncut bubble profile on the whole manifold. The formulas in
the core, including the derivatives, converge after rescaling
to the corresponding Euclidean formulas. The same bounds hold
after multiplication by $\theta_y$. The positivity and two-sided
profile bound for $V$ imply the bounds by $C U_\epsilon^p$ for
the derivatives of the power, including the term with $V^{p-2}$
in a second derivative. Here the profiles in the matching and
exterior regions use distance on the manifold.

The critical tails of these bounds are uniformly small. Thus
differentiating the critical mass once and twice and rescaling
gives the corresponding Euclidean derivatives, which vanish
by invariance of the full-space mass under translation and
dilation. This gives convergence of the normalized source modes
in $L^s$. Applying $L_g^{-1}$ gives convergence of the projected
modes in $H^1$ to the cutoff Euclidean modes. To see the latter
comparison, the Laplacian of a cutoff Euclidean mode agrees
with its differentiated Euclidean source up to an $o(1)$ error
in $H^{-1}$. On a ball of fixed small radius the principal
coefficient error is $O(r^2)$ times a bounded Dirichlet norm,
and outside the ball that norm tends to zero. Lower-order and
cutoff errors tend to zero by the $L^2$ bound and the Dirichlet
tail. Taking the radius to zero after the scale limit proves the
comparison at both derivative orders.
On the unit round sphere use the exact family
$$
V^{\rm rd}_{\epsilon,y}(x)=K
\left(\frac{\epsilon}{\epsilon^2\cos^2(d(x,y)/2)
+4\sin^2(d(x,y)/2)}\right)^a
$$
Stereographic radius $r=2\tan(d/2)$ proves both
$\int(V^{\rm rd})^{p+1}=1$ and
$L_{\rm rd}V^{\rm rd}=Y(V^{\rm rd})^p$ by conformal covariance.
If $g=\varphi^{4/(n-2)}F^*g_{\rm rd}$, use the round scale
$\epsilon_{\rm rd}=\varphi(y)^{-2/(n-2)}\epsilon$ and transfer
by $\varphi^{-1}F^*$. The rescaled core on $g$ has exactly scale
$\epsilon$, since the conformal distance factor at the center is
$\varphi(y)^{-2/(n-2)}$. Take $f=V^p$ and $B^g=V$.
Then $D_g=0$ and the residual vanishes identically. Differentiating
this explicit positive global family proves its stated mode bounds;
there is no truncation error at zero mass.
\end{proof}

\subsection{Comparison of the Green projections}

\begin{lemma}\label{all-potential-interactions}
Under $Y_\Psi(M)>0$ and $q=R_g-\Psi\ge0$, $q\not\equiv0$, use exactly
the sources in Proposition \ref{all-positive-source}, including the
exact family in the round case, and set $B=YL_\Psi^{-1}f$.
They satisfy all the hypotheses of
Proposition \ref{source-convergence}, for $n\ge4$.
\end{lemma}

\begin{proof}
The resolvent identity gives
$$
D_\Psi=D_g+Y^{-2}S,\qquad
S=\int_Mq B^gB,\qquad B^g\le B\le CB^g.
$$
The positive minimum of each global Green kernel and
$\int f\ge c\epsilon^a$ give
$S\ge c\epsilon^{2a}\int_Mq$ for every center.
For the equation error split the source $qB^g$ into a ball of fixed
radius $2\rho$ and its complement. The near inverse has $H^1$ norm
at most
$C\|\sqrt q\|_{L^n(B_{2\rho})}\|\sqrt q B^g\|_2\le C\rho\sqrt S$.
Away from the ball $qB^g$ is bounded by $C_\rho\epsilon^a$.
Since $\sup_x\int_M G_\Psi(x,z)\,dv_g(z)<\infty$, the far
inverse has $L^\infty$ norm at most $C_\rho\epsilon^a$.
Its contribution to the power difference is at most
$$
C_\rho\epsilon^a\|(B^g)^{p-1}\|_s
+C_\rho\epsilon^{a+2}=o_\rho(\epsilon^a).
$$
Here the power norm tends to zero by its subcritical exponent and the
uniform critical profile bound. Square, take the scale limit, and then
let $\rho\downarrow0$. Proposition \ref{all-positive-source} and the
global lower bound for $S$ prove
$\|L_\Psi B-YB^p\|_s^2=o(D_\Psi)$, uniformly in the center, including
centers approaching the zero set at any speed. For the limit of the gain, the profile bound gives
$$
 \|B^g\|_2^2\le C
 \begin{cases}
 \epsilon^2(1+|\log\epsilon|),&n=4,\\
 \epsilon^2,&n\ge5.
 \end{cases}
$$
Indeed the exterior radial integral is bounded by
$C\epsilon^{n-2}\int_\epsilon^{\operatorname{diam}M}r^{3-n}\,dr$
and the core by $C\epsilon^2$.
$S\le C\|B^g\|_2^2\to0$ and $D_\Psi\to0$ in every dimension.
The compact operator
$L_\Psi^{-1}q:H^1\to H^1$ takes bounded weakly vanishing geometric
parameter modes to $o(1)$, proving the parameter assertions.
The source profile and its lower bound on $B_\epsilon(y)$
follow from the same construction: a fixed smaller ball in the
rescaled conformal coordinates lies in $B_\epsilon^g(y)$ and
carries at least $c\epsilon^a$ source mass.
This verifies \eqref{family-profile}--\eqref{family-error} and the
two parameter orders of Proposition \ref{source-convergence}.

\end{proof}

\begin{proof}[Proof of the convergence implication in Theorem \ref{main-iff}]
For a nonround conformal class and $\Psi=R_g$,
Proposition \ref{all-positive-source} gives the source and all
the projection estimates in Proposition \ref{source-convergence}.
The pointwise profile and the local source mass in
\eqref{family-profile} follow from the rescaled core and the Green
tail in its construction. If $\Psi\le R_g$ and $\Psi\not\equiv R_g$,
Lemma \ref{all-potential-interactions} gives the hypotheses of the
same criterion, including the case of a round conformal class.

It remains to take $\Psi=R_g$ in the round conformal class.
With geometric time $4t/(n-2)$ and initial metric
$u_0^{4/(n-2)}g$, \eqref{flow} is the volume-normalized Yamabe
flow. Ye's convergence theorem in the locally conformally flat
case \cite{Ye94} applies. In each case the limit satisfies
\eqref{stationary}.
\end{proof}

\section{Convergence in dimension three}\label{projective-section}

The regular part \eqref{three-mass} enters the elliptic compactness
theory in dimension three \cite{LZ99}. Here its sign controls the
energy of the Green projections in Proposition~\ref{source-convergence}.
In a fixed small normal ball about $y$, consider
$\chi(d_g(x,y))/(32\pi d_g(x,y))$, where $\chi$ is one near zero
and vanishes outside the ball. Applying $L_\Psi$ gives
$\delta_y$ plus a function bounded by $C/d_g(x,y)$. This follows
from the radial Laplacian formula
$\Delta_g d_g(\cdot,y)=2/d_g(\cdot,y)+O(d_g(\cdot,y))$.
Applying $L_\Psi^{-1}$ to the other term shows that
$$
 G_\Psi(x,y)=\frac{1}{32\pi\,d_g(x,y)}+H_\Psi(x,y)
$$
near the diagonal, with $H_\Psi$ continuous there. Indeed the
convolution of two kernels bounded by inverse distances is
jointly continuous. On the union of radius $\rho$ balls about
the two poles its integral is $O(\rho)$ by Cauchy--Schwarz
and $\int_{B_\rho(z)}d_g(x,w)^{-2}\,dv_g(w)\le C\rho$,
uniformly in $z,x$. On the complement of those balls the
integrands are smooth. This also proves continuity of
$A_\Psi(y)=H_\Psi(y,y)$.

\begin{proof}[Proof of Theorem \ref{main-three}]
Suppose first that $A_\Psi>0$. Fix a cutoff supported in a
small normal ball, equal to one on a smaller ball.
In normal coordinates at $y$ let $b_{\epsilon,y}$ be the
cutoff of $U_\epsilon$ normalized to have $L^6(M)$ norm one.
Set $f=b_{\epsilon,y}^5$ and $B=YL_\Psi^{-1}f$. The
source bounds in \eqref{family-profile} follow by rescaling.
The $L^6$ normalizing factor is $1+O(\epsilon^2)$, by the
$O(|x|^2)$ volume error and the $O(\epsilon^3)$ Euclidean
critical mass tail.

We compute the difference in \eqref{family-gain}.
For $z,w$ in these coordinates the difference between
$d_g(\exp_y z,\exp_y w)^{-1}$ and $|z-w|^{-1}$ is bounded
by $C(|z|^2+|w|^2)|z-w|^{-1}$. To see this, compare lengths
in the ball of radius $3\max(|z|,|w|)$, which contains the
minimizing segment, using the second-order metric error.
The straight coordinate segment gives the reverse length comparison.
The Newton kernel applied to $U_\epsilon^5$ gives
$U_\epsilon/Y$. Consequently the errors from the distance
and volume comparisons in the double integral are at most
$C\int_{\mathbb R^3}|z|^2 U_\epsilon^6\,dz=O(\epsilon^2)$.
Normalization has the same error. The deleted Newtonian
interaction with the exterior source is bounded by
$2Y^{-1}\int_{\{\chi\ne1\}}U_\epsilon^6=O(\epsilon^3)$,
by positivity. The uncut Newtonian double integral is $1/Y$.
For the continuous part of the kernel use
$$
 \epsilon^{-1/2} f\,dv_g\rightharpoonup
 \left(\int_{\mathbb R^3} U_1^5\,dz\right)\delta_y
$$
in centered coordinates. Uniform continuity of the regular
part near the diagonal and the uniform source mass estimates
then give
$$
 D_\Psi(\epsilon,y)=\epsilon A_\Psi(y)
       \left(\int_{\mathbb R^3} U_1^5\,dz\right)^2+o(\epsilon)
$$
uniformly in $y$. In particular $D_\Psi\ge c\epsilon$.

On a smaller core ball the same kernel comparison yields
$B=U_\epsilon+O(\sqrt{\epsilon})$ uniformly. The contribution
from the regular part is bounded by $C\int_M f$. For the
Newtonian distance error, the term with $|z|^2$ is at most
$C|z|^2U_\epsilon(z)\le C\sqrt{\epsilon}$. The term with
$|w|^2$ is bounded by the Newton potential of
$C\epsilon^{5/2}(\epsilon^2+|w|^2)^{-3/2}$, which is radial
decreasing and has value $O(\epsilon^{3/2})$ at the origin.
The cutoff and normalization errors are no larger. On the
complement of the smaller core, the general Green upper bound
gives $B=O(\sqrt{\epsilon})$. Hence
$$
 \|L_\Psi B-YB^5\|_{6/5}
       =Y\|f-B^5\|_{6/5}\le C\epsilon .
$$
On the core the power difference is at most
$C(\epsilon^2U_\epsilon^5+\sqrt{\epsilon}U_\epsilon^4+
\epsilon^{5/2})$. Its $L^{6/5}$ norm is $O(\epsilon)$ because
$\|U_\epsilon^4\|_{6/5}=O(\sqrt{\epsilon})$.
All terms outside the core are $O(\epsilon^{5/2})$. This is
the squared-residual comparison in \eqref{family-error}.

For the parameter estimates, write $b_{\epsilon,y}$ in terms
of the squared distance in the normal ball. For any smooth
radial expression $F(r^2)$ there,
$$
 \Delta_g F(r^2)-\Delta_{\mathbb R^3}F(r^2)
       =(\Delta_g r^2-6)F'(r^2).
$$
The smooth coefficient $\Delta_g r^2-6$ vanishes to second
order in the normal variable. Each normalized center
derivative costs at most $C\epsilon/(\epsilon+r)$ on the
profile. Differentiating the vanishing coefficient once or
twice still bounds the product by $CU_\epsilon$ after the
normalizing factors $\epsilon$ are included. Scale derivatives
of the rational profile have the same order. The normalizing
factor and its first two derivatives differ from those of the
Euclidean mass by $O(\epsilon^2)$, from the same volume
integral after rescaling. On the cutoff annulus all derivatives
under consideration have size $O(\sqrt\epsilon)$. Thus
$L_\Psi b_{\epsilon,y}-Yf$ and its first two normalized
parameter derivatives have $L^{6/5}$ norm $O(\sqrt\epsilon)$.
Applying $L_\Psi^{-1}$ proves convergence of the projected
modes in $H^1$ to the cutoff Euclidean ones. Proposition
\ref{source-convergence} applies.

For the assertion when $\Psi\le R_g$, the conformal Laplacian
is coercive because $Y(M,[g])\ge Y_\Psi(M)>0$.
The geometric regular part $A_{R_g}$ is nonnegative, and it
is positive if the conformal class is not round. In
conformal normal coordinates the expansion of the normalized
Green function is $r^{-1}+32\pi A_{R_g}(y)+O''(r)$
\cite[Lemma 6.4]{LP87}. The conformal factor is one with zero
first derivative at the pole, so the change of distance and
the Green transformation preserve the regular part. The
scalar-flat Green metric has ADM mass $64\pi A_{R_g}(y)$.
In the orientable case its order-one decay through two derivatives
allows any weighted Sobolev decay rate in $(1/2,1)$ for the
time-symmetric application of \cite[Theorem~1 and Section~2]{EHLS16}.
Rigidity follows from \cite[Lemma~10.7]{LP87} as in
Lemma~\ref{all-mass-flux}. For a nonorientable manifold the
single-pole comparison on the orientable double cover again applies.

The resolvent identity, evaluated on the diagonal after
subtracting the fundamental singularity, gives
$$
 A_\Psi(y)-A_{R_g}(y)
 =\int_M G_\Psi(y,z)(R_g-\Psi)(z)G_{R_g}(z,y)\,dv_g(z).
$$
The integral is finite and continuous by the same small-ball
estimate used above. If $R_g-\Psi\not\equiv0$, the integral
is strictly positive at every $y$. Thus the first part of
the proof applies both in that case and when the conformal
class is nonround. In the remaining case $\Psi=R_g$ in the
round class, use Ye's theorem \cite{Ye94} as above.
\end{proof}

\begin{proof}[Proof of Corollary \ref{main-counterexample}]
The potential is positive, so the operator is coercive and
$Y_\Psi(M)>0$. Put $\sigma=\sqrt{1-\Psi/8}$ and
$\rho=d_g(x,y)$. The positive Green kernel on unit round
$\mathbb{RP}^3$ is
$$
 G_\Psi(x,y)=\frac{1}{32\pi}
   \frac{\cos(\sigma(\rho-\pi/2))}
        {\cos(\pi\sigma/2)\sin\rho},\qquad 0<\rho\le\pi/2.
$$
Indeed, writing a radial function as $v(\rho)/\sin\rho$
reduces the homogeneous equation to $v''+\sigma^2v=0$.
The numerator is even about $\rho=\pi/2$. On the spherical
cover it can be written near the cut locus as
$\cos(\sigma\arcsin\langle x,y\rangle)$, and the denominator
there is a positive multiple of
$\sqrt{1-\langle x,y\rangle^2}$. Thus the quotient is smooth
across the cut locus. The singular coefficient gives unit
flux at the single pole. Expansion at that pole yields
$$
 A_\Psi=\frac{\sigma\tan(\pi\sigma/2)}{32\pi}>0 .
$$
Theorem \ref{main-three} proves convergence. On this metric
$R_g=6$, so the pointwise condition in Theorem \ref{main-iff}
does not characterize convergence in dimension three.
\end{proof}

\section{Concentration at a prescribed point}\label{four-section}

Suppose $n\ge4$. Fix a compact set of centers in a coordinate
chart. Take a fixed radial cutoff $\chi$ equal to one near zero
and supported within the injectivity radius. With the profile
\eqref{euclidean-bubble}, define
$$
 b_{\epsilon,y}
 =\frac{\chi(d_g(\cdot,y))U_\epsilon(d_g(\cdot,y))}
        {\|\chi(d_g(\cdot,y))U_\epsilon(d_g(\cdot,y))\|_{p+1}} .
$$
The radial expression is smooth in $x$ and $y$, since the
squared distance is smooth near the diagonal and the cutoff
vanishes before the cut locus. Write $D=\epsilon\partial_\epsilon$,
and set
$$
 \Phi(\epsilon)=
 \begin{cases}
 \epsilon^2\log(1/\epsilon),& n=4,\\
 \epsilon^2,& n\ge5.
 \end{cases}
$$

\begin{lemma}\label{four-bubble-estimates}
Uniformly in the center, for $C_4=2\pi\sqrt6$ and
$C_n=\int_{\mathbb R^n}U_1^2\,dx$ if $n\ge5$, one has
\begin{align}\label{four-bubble-energy}
 E_\Psi(b_{\epsilon,y})-Y
   &=C_n(\Psi-R_g)(y)\Phi(\epsilon)+o(\Phi(\epsilon)),\notag\\
 D E_\Psi(b_{\epsilon,y})
   &=2C_n(\Psi-R_g)(y)\Phi(\epsilon)+o(\Phi(\epsilon)).
\end{align}
For $R_b=L_\Psi b_{\epsilon,y}-Yb_{\epsilon,y}^p$,
\begin{equation}\label{four-bubble-residual}
 \|R_b\|_s+\|D R_b\|_s\le C\omega(\epsilon),\qquad
 \omega(\epsilon)=
 \begin{cases}
 \epsilon,&n=4,\\
 \epsilon^{3/2},&n=5,\\
 \epsilon^2(1+|\log\epsilon|)^{2/3},&n=6,\\
 \epsilon^2,&n\ge7.
 \end{cases}
\end{equation}
In particular $\omega^2=o(\Phi)$. The profile satisfies
$|Db_{\epsilon,y}|\le Cb_{\epsilon,y}$ and
$\int d_g(x,y)^2 b_{\epsilon,y}^{p+1}\,dv_g\le C\epsilon^2$.
Its normalized scale and center derivatives through order two
are bounded in $H^1$ and converge, after rescaling at the center,
to the Euclidean modes with vanishing Dirichlet tails.
\end{lemma}

\begin{proof}
Write $r=|x|$ in normal coordinates about $y$. The radial
Dirichlet density is exactly $|\partial_r U_\epsilon|^2$.
The angular average of the volume density is
$1-R_g(y)r^2/(6n)+O(r^3)$. The cutoff changes the Euclidean
Dirichlet integral by $O(\epsilon^{n-2})$ and the critical
mass by $O(\epsilon^n)$. The second moment of $U_1^{p+1}$
is finite, so the mass-normalizing factor is
$1+O(\epsilon^2)$.

For $n\ge5$ the potential contribution has leading term
$\Psi(y)\epsilon^2\int U_1^2$. To compute the geometric
term, test $-c_n\Delta U_1=YU_1^p$ with $r^2U_1$ and
integrate the radial derivative in the cross term. This gives
$$
 c_n\int r^2|\nabla U_1|^2
    =Y\int r^2 U_1^{p+1}+n c_n\int U_1^2,
$$
with all integrals on $\mathbb R^n$. The boundary term
tends to zero for $n>4$. Direct integration of the
two radial powers gives
$$
 Y\int r^2 U_1^{p+1}
   =\frac{n^2(n-4)}{n-2}\int U_1^2.
$$
Indeed the ratio of the two radial integrals for
$r^2 U_1^{p+1}$ and $U_1^2$, before the bubble constants,
is $(n/2)(n/2-2)/[(n-1)(n-2)]$, and
$YK^{p-1}=4n(n-1)$ for the constant $K$ in $U_1$.
Combining the Dirichlet and normalization terms gives
$$
 c_n\int r^2|\nabla U_1|^2
   -\frac{2Y}{p+1}\int r^2U_1^{p+1}=6n\int U_1^2.
$$
Multiplication by $-R_g(y)/(6n)$ gives the remaining
coefficient of $\epsilon^2$. The same leading energy
term in dimensions $n\ge6$, in the action normalization,
is computed in \cite[Lemma 3.1]{ChenMa26}.

When $n=4$, $U_\epsilon=K\epsilon/(\epsilon^2+r^2)$ and
$|S^3|K^2=2\pi\sqrt6$. The potential and quadratic volume
corrections give respectively
\begin{align*}
 &|S^3|K^2\Psi(y)\epsilon^2
       \int_0^\delta\frac{r^3}{(\epsilon^2+r^2)^2}\,dr,\\
 &-|S^3|K^2 R_g(y)\epsilon^2
       \int_0^\delta\frac{r^7}{(\epsilon^2+r^2)^4}\,dr,
\end{align*}
up to $O(\epsilon^2)$, where the cutoff is one on $r\le\delta$.
Both radial integrals are $\log(1/\epsilon)+O(1)$ by
splitting at $r=\epsilon$. The higher volume terms, the
variation of the potential, and the annular terms have size
$O(\epsilon^2)$, also after applying $D$. The mass factor
therefore changes only this error. Each $D$ derivative
of the two displayed integrals is bounded, as follows by
differentiating under the integral. This gives the scale
derivative in dimension four. For $n\ge5$, rescale
$r=\epsilon |z|$ in the volume and potential integrals
before differentiating. After division by $\epsilon^2$,
the quadratic metric terms and their derivatives are
dominated by the integrable second moments of
$|\nabla U_1|^2$ and $U_1^{p+1}$. The potential terms
are dominated by $C U_1^2$. This proves the differentiated
higher-dimensional expansion by dominated convergence.

In the core,
$|(\Delta_g-\Delta)U_\epsilon|\le Cr|\partial_r U_\epsilon|
\le C U_\epsilon$. This estimate also holds with $DU_\epsilon$,
and the potential term satisfies the same bound. On the
fixed cutoff annulus all equation errors have size
$O(\epsilon^a)$. Normalization contributes $O(\epsilon^2)$
in $L^s$. Integrating $U_\epsilon^s$ over the fixed core
gives \eqref{four-bubble-residual}, with the logarithm
when $n=6$. The pointwise scale bound and moment bound
follow from the radial formula and its critical second moment.

For the parameter derivatives, the first center derivative
of $r^2$ is $O(r)$ and the second is bounded. The chain
rule bounds the normalized derivatives through order two
by $C U_\epsilon$, and their spatial gradients by
$C U_\epsilon/(\epsilon+r)$, including cutoff terms.
With $r>A\epsilon$ the squared gradient integrals are at
most $C A^{2-n}$. On a fixed rescaled ball the metric tends
to the Euclidean metric, the mass factor tends to one
with its derivatives tending to zero, and the rational
profile and its modes converge there. This proves the
parameter assertions.
\end{proof}

Let the center range over a compact coordinate set on which
$\Psi-R_g>0$. A function in a uniform $H^1$ neighborhood
of $b_{\epsilon,y}$ can be written
$$
 u=\alpha b_{\widetilde\epsilon,\widetilde y}+w
$$
with the remainder orthogonal in $\langle\cdot,\cdot\rangle_\Psi$
to amplitude, scale and translation modes. In fact, in
the coordinates $(\log(\widetilde\epsilon/\epsilon),
(\widetilde y-y)/\epsilon)$ the Jacobian of the orthogonality
equations is the negative Gram matrix of the modes at
$w=0$ and $\alpha=1$. The last assertion of Lemma
\ref{four-bubble-estimates} gives a uniformly invertible
Euclidean limit and uniform second parameter bounds.
The implicit function theorem yields a local branch with
$\alpha$ close to one, a small relative change of scale,
and a center displacement small relative to $\epsilon$.
At any later representation on this branch, the last
$n+1$ columns of the leading Gram matrix are multiplied
by $\alpha$, and the remaining Jacobian terms are
$O(\|w\|_\Psi)$ by the second derivative bounds. Hence
the branch can be continued whenever these quantities
stay within the given neighborhood. In particular at a
finite endpoint it suffices that the scale be bounded below
and that the upper scale bound, center region, remainder
bound and amplitude range retain a strict margin. Under
these conditions, ordinary parameter derivatives are bounded
and differentiation of the orthogonality equations makes
the parameter speeds bounded by $C\|u_t\|_{H^1}$, so the
parameters have a limit at the endpoint.

On the orthogonal complement, uniformly for small scales,
\begin{equation}\label{four-coercive}
 \|w\|_\Psi^2-pY\int b_{\epsilon,y}^{p-1}w^2
 \ge c\|w\|_\Psi^2 .
\end{equation}
Indeed, for a sequence with $\|w\|_\Psi=1$, rescale at
the center and take a weak Dirichlet limit. Its energy is
at most one, by weak lower semicontinuity after subtracting
the global weak limit. The global limit contributes zero
to the weighted integral by local approximation with
smooth functions and $\|b^{p-1}\|_{L^1}\to0$.
The tails of $b^{p-1}$ are uniformly small in the rescaled
$L^{n/2}$ norm and Rellich compactness applies on rescaled
balls. Orthogonality passes to the Euclidean modes by the
Dirichlet tail estimates. Equation \eqref{euclidean-gap}
then bounds the negative term by $1-4/(n+4)$ in the
limit, which proves \eqref{four-coercive}.

\begin{proposition}\label{four-static}
For centers in a compact subset of $\{\Psi>R_g\}$, take
$u=\alpha b_{\epsilon,y}+w\ge0$ of critical mass one with
the preceding orthogonality conditions and $1/2\le\alpha\le2$.
For all sufficiently small $\epsilon$ and $\|w\|_\Psi$,
$$
 E_\Psi(u)-Y\asymp \Phi(\epsilon)+\|w\|_\Psi^2 .
$$
For smooth $u$ one also has
\begin{equation}\label{four-residual-lower}
 E_\Psi(u)-Y\le C\|L_\Psi u-E_\Psi(u)u^p\|_s .
\end{equation}
\end{proposition}

\begin{proof}
Write $b=b_{\epsilon,y}$. Since $w$ is orthogonal to $b$,
$\int b^pw=-Y^{-1}\int R_bw=O(\omega\|w\|_\Psi)$.
Taylor expansion of the mass first gives
$|\alpha-1|\le C(\omega\|w\|_\Psi+\|w\|_\Psi^2)$,
and then gives
$$
 \alpha^2=1-2\int b^pw-p\int b^{p-1}w^2
              +o(\|w\|_\Psi^2).
$$
The last remainder is uniform as the scale and the
remainder norm tend to zero. For $p+1<3$ the power
remainder in the mass is bounded by
$C\|w\|_{p+1}^{p+1}$, and for $p+1\ge3$ it is bounded
by $C(\|w\|_{p+1}^3+\|w\|_{p+1}^{p+1})$. The square of
the first-order amplitude error is of order
$\omega^2\|w\|_\Psi^2+\|w\|_\Psi^4$ and is included in
the remainder. Since $E_\Psi(u)=\alpha^2E_\Psi(b)+\|w\|_\Psi^2$,
\eqref{four-coercive} controls the transverse quadratic
part. Young's inequality absorbs the linear error at a
cost $C\omega^2=o(\Phi)$, which is smaller than the
positive term in \eqref{four-bubble-energy}. This proves
the two-sided energy estimate.

Denote the norm on the right of \eqref{four-residual-lower}
by $R$. Pair the residual with $w$. Its quadratic part
is $\|w\|_\Psi^2-pE_\Psi(u)\alpha^{p-1}\int b^{p-1}w^2$,
which has a positive gap by \eqref{four-coercive}, since
$E_\Psi(u)\alpha^{p-1}\to Y$. The contribution of
$b^p$ is $O(\omega\|w\|_\Psi)$. The remaining terms after
pairing are $o(\|w\|_\Psi^2)$ by the power Taylor bound
and H\"older. Consequently
$$
 \|w\|_\Psi\le C(R+\omega).
$$
For the scale pairing use normalization $\int b^p Db=0$
and the identity
$L_\Psi Db-pY b^{p-1}Db=D R_b$. Orthogonality to $Db$
then gives
\begin{align*}
 \int (L_\Psi u-E_\Psi(u)u^p)Db
 ={}&\frac{\alpha}{2}D E_\Psi(b)
 +\frac{E_\Psi(u)\alpha^{p-1}}{Y}\int (D R_b)w
 +O(\|w\|_\Psi^2).
\end{align*}
For $1<p<2$ the power remainder where $b>0$ is at
most $C b^{p-2}|w|^2$, since $\alpha b+w\ge0$ and
$\alpha$ lies in a fixed positive interval. Pairing
with $|Db|\le Cb$ bounds it by $C\int b^{p-1}w^2$.
On $\{b=0\}$ the pairing vanishes. For $p\ge2$ one
also has a term bounded by $C\int b|w|^p$, of order
at most $\|w\|_\Psi^2$ for small remainder. Thus
\eqref{four-bubble-energy} gives
$\Phi\le C(R+R^2+\omega^2)$, using the preceding bound
on $w$. Absorb $\omega^2=o(\Phi)$. For $R\le1$ the
two-sided energy estimate now gives \eqref{four-residual-lower}.
For $R>1$ it follows by shrinking the neighborhood so the
energy difference is less than one.
\end{proof}

\begin{proof}[Proof of concentration in Theorem \ref{main-iff}]
Choose normal coordinates about the prescribed point $y$,
and fix a closed coordinate ball $\mathcal K$ about $y$
in $\{\Psi>R_g\}$. Choose a smooth map $F:M\to\mathbb R^n$
which agrees with these coordinates on a neighborhood of
$\mathcal K$, with $F(y)=0$. Take the initial bubble center
$z$ in a smaller closed coordinate ball about $y$. For
$\sigma>0$ set
$$
 u_0=\frac{b_{\sigma,z}+\sigma^4}{\|b_{\sigma,z}+\sigma^4\|_{p+1}} .
$$
Its distance from $b_{\sigma,z}$ in $H^1$ is $O(\sigma^4)$.
By \eqref{four-bubble-energy} its energy is above $Y$ and
tends to $Y$ uniformly in $z$. Take $\sigma$ small enough
that the energy excess is less than $\eta$ and continue
the local modulation for the global solution from
Proposition \ref{evolution}.

As long as the modulation remains in the region of
Proposition \ref{four-static}, write $e=r-Y$ and
$\mathcal D=\int u^{p-1}u_t^2$. One has
$$
 0<e\asymp\Phi(\epsilon)+\|w\|_\Psi^2,\qquad
 e\le C\sqrt{\mathcal D},\qquad e'=-2\mathcal D.
$$
The middle bound follows from \eqref{four-residual-lower}
and $\|L_\Psi u-r u^p\|_s\le\sqrt{\mathcal D}$ by mass-one
H\"older. The energy bound keeps the scale below the upper
scale threshold, because $\Phi$ is increasing for small
$\epsilon$, and it keeps the remainder small. The mass
constraint then bounds $\alpha$ in a fixed compact interval
near one.

For the center put $m(t)=\int_M F u^{p+1}$. By
Lemma \ref{four-bubble-estimates}, Sobolev's inequality,
and the energy estimate,
$$
 \int |F-F(y(t))|^2u^{p+1}
     \le C(\epsilon^2+\|w\|_\Psi^{p+1})\le C e ,
$$
where $y(t)$ denotes the modulated center. Hence
$|m(t)-F(y(t))|\le C\sqrt e$. By mass preservation,
$m'=(p+1)\int (F-m)u^p u_t$. The mean minimizes the
variance, so Cauchy--Schwarz gives
$$
 |m'|\le C\sqrt e\sqrt{\mathcal D}
          \le C\frac{\mathcal D}{\sqrt e}
          =-C\frac{d}{dt}\sqrt e.
$$
Thus $m(t)$ moves by at most $C\sqrt{e(0)}$ throughout
the continuation interval. The initial bubble moment
gives $|m(0)-F(z)|\le C\sigma$ by the same variance
calculation with the additive perturbation. These bounds
keep the center strictly inside $\mathcal K$ for every
$z$ in the smaller ball when $\sigma$ is sufficiently
small. All parameter boundaries except $\epsilon=0$
therefore retain strict margins.

At any finite time the global solution is bounded.
Uniformly in the center, a sufficiently large fixed
rescaled ball $B_{A\epsilon}(y(t))$ carries a positive
amount of the mass of $b_{\epsilon,y(t)}$. With $\alpha$
near one and $\|w\|_{p+1}$ small, it carries a fixed
positive amount of the mass of $u$. This excludes any
sequence of scales tending to zero at a finite endpoint,
using the finite-time bound on $u$. The local modulation
thus continues for all time.

Now $e'\le-ce^2$ gives $e(t)\le C/(1+t)$, uniformly for
the chosen family of initial centers. Consequently
$\epsilon(t)\to0$, $\|w(t)\|_\Psi\to0$ and $\alpha(t)\to1$.
The barycenter has finite total variation with tail
bounded by $C\sqrt{e(t)}$. Its limit lies in the coordinate
image of $\mathcal K$ and is also the limit of $F(y(t))$.
Thus $y(t)\to y_\infty$ for a point $y_\infty$ in
$\mathcal K$. The bubble moment and vanishing remainder
imply $u(t)^{p+1}\,dv_g\rightharpoonup\delta_{y_\infty}$.
The mass in $B_{A\epsilon}(y(t))$ and its volume
$O(\epsilon^n)$ give $\|u(t)\|_\infty\ge c\epsilon^{-a}$.
In dimension four the lower bound for $e$ gives
$\epsilon(t)^2\log(1/\epsilon(t))\le C/(1+t)$.
For large $t$ this first gives $\epsilon(t)\le Ct^{-1/2}$ and then
$\log(1/\epsilon(t))\ge c\log t$. As $a=1$ in this dimension,
the resulting bounds are
$$
\epsilon(t)\le C(t\log t)^{-1/2},\qquad
\|u(t)\|_\infty\ge c(t\log t)^{1/2}.
$$

To prescribe $y_\infty$, keep the same $\sigma$ for all
initial centers in the smaller ball. The limit $m(\infty)$
depends continuously on the initial center: finite-time
solutions depend continuously on the data and the tail of
$m$ is bounded uniformly by $C(1+t)^{-1/2}$. The map
$F(z)\mapsto F(y_\infty(z))$ differs from the identity
by at most $C(\sigma+\sup_z\sqrt{e_z(0)})$. Reduce
$\sigma$ so that this error is less than the radius of
the smaller coordinate ball. Then $F(z)\mapsto
F(z)-F(y_\infty(z))$ maps the closed ball into itself.
Brouwer's theorem gives a fixed point, for which
$F(y_\infty(z))=0$. This yields concentration at the
prescribed point $y$ and completes the proof of
Theorem \ref{main-iff}.
\end{proof}

\medskip
\noindent\textbf{Data Availability Statement.}
Data sharing is not applicable to this article because no datasets
were generated or analyzed during the current study.


\noindent\textbf{Conflict of interest.}
The authors declare that they have no conflict of interest.

\noindent\textbf{Code availability.}
Not applicable.

\noindent\textbf{AI Assistance Statement.}
AI is used as an assistive tool.

\end{document}